\documentclass[11pt,letterpaper]{amsart}

\usepackage[utf8]{inputenc}
\usepackage[T1]{fontenc}
\usepackage{lmodern}
\usepackage[english]{babel}
\usepackage{amsmath,amssymb,amsthm,mathtools}
\usepackage[stretch=50,shrink=50]{microtype}
\usepackage{enumitem}
\usepackage{needspace}
\usepackage{xurl}
\usepackage[colorlinks,citecolor=blue,linkcolor=blue]{hyperref}

\numberwithin{equation}{section}
\newtheorem{theorem}{Theorem}[section]
\newtheorem{proposition}[theorem]{Proposition}
\newtheorem{lemma}[theorem]{Lemma}
\newtheorem{corollary}[theorem]{Corollary}

\theoremstyle{definition}
\newtheorem{definition}[theorem]{Definition}
\newtheorem{condition}{Condition}

\theoremstyle{remark}
\newtheorem*{remark}{Remark}

\newcommand{\R}{\mathbb R}
\newcommand{\Q}{\mathbb Q}
\newcommand{\Z}{\mathbb Z}

\newcommand{\N}{\mathbb N}

\newcommand{\gbar}{\bar g}

\DeclareMathOperator{\Ric}{Ric}
\DeclareMathOperator{\Scal}{Scal}
\DeclareMathOperator{\Hess}{Hess}
\DeclareMathOperator{\Vol}{Vol}
\DeclareMathOperator{\tr}{tr}
\DeclareMathOperator{\Hom}{Hom}

\DeclareMathOperator{\supp}{supp}
\DeclareMathOperator{\ind}{ind}
\DeclareMathOperator{\Crit}{Crit}
\DeclareMathOperator{\intt}{int}
\DeclareMathOperator{\AVR}{AVR}
\DeclareMathOperator{\Ext}{Ext}

\title[RICCI CURVATURE AND CONTRACTIBLE 4-MANIFOLDS]%
{Contractible 4-Manifolds\\ with Nonnegative Ricci Curvature\\ and Euclidean Volume Growth}

\author{Junhyun An}
\address{Seoul National University, Seoul, Republic of Korea}
\email{anjunhyun@snu.ac.kr}

\author{Seunghwan Yoon}
\address{Korea Advanced Institute of Science and Technology,
Daejeon, Republic of Korea}
\address{RTst, Daejeon, Republic of Korea}
\email{seunghwanyoon@kaist.ac.kr}


\hypersetup{
  pdftitle={Contractible 4-Manifolds with Nonnegative Ricci Curvature
    and Euclidean Volume Growth},
  pdfauthor={Junhyun An, Seunghwan Yoon},
  pdfsubject={Green potentials, PSC cuts, simply connected exteriors}
}

\begin{document}

\begin{abstract}
We prove that every contractible smooth 4-manifold admitting
a complete Riemannian metric with nonnegative Ricci curvature
and Euclidean volume growth is homeomorphic to $\R^4$.
More precisely, every compact subset lies in the interior of
a connected compact regular domain with simply connected
complement.
The key construction produces compact domains whose boundaries
admit positive scalar curvature and which are regular sublevels
of global strictly three-convex Morse exhaustions.
A topological criterion forces the exterior fundamental group
to be trivial whenever it contains no nonabelian free subgroup
and satisfies the resulting boundary, homological, and relative
cell hypotheses; in the geometric application, this group is finite.
Its proof uses a cellular formulation of the
Gerstenhaber--Rothaus theorem and a smooth filling obstruction.
The resulting complements give simple connectivity at infinity,
and Freedman's recognition theorem completes the argument.
\end{abstract}
\maketitle

\vspace{1.5em}

\section{Introduction}\label{sec:intro}

Contractibility alone does not determine the topology of an open
4-manifold. Although every loop contracts in the ambient manifold,
loops far from a compact set need not contract while avoiding that
set. Controlling this behavior at infinity is therefore central to
recognizing a contractible 4-manifold as $\R^4$.
In this context, Semola conjectures that every contractible
4-manifold admitting a complete metric with nonnegative Ricci
curvature and Euclidean volume growth is homeomorphic to
$\R^4$ \cite[Conjecture~5.8]{Semola}.
We prove this conjecture.

\begin{condition}[Geometric hypotheses]\label{cond:geometry}
Let $M$ be a smooth 4-manifold without boundary, equipped with a
smooth Riemannian metric $g$. Assume that $M$ is contractible,
$(M,g)$ is complete, and $\Ric_g\ge0$. Fix $p\in M$ and assume
\[
 \AVR(g):=\lim_{R\to\infty}
 \frac{\Vol_gB_R^g(p)}{\omega_4R^4}>0,
\]
where $\omega_4$ is the volume of the Euclidean unit ball in $\R^4$.
In particular, $M$ is connected and noncompact.
\end{condition}

\begin{theorem}\label{thm:main}
Under Condition~\ref{cond:geometry}, $M\cong_{\mathrm{top}}\R^4$.
\end{theorem}

More precisely, we obtain the following control of compact
complements, writing $O_\Gamma=M\setminus\Gamma$.

\begin{theorem}[Cofinal simply connected complements]
\label{thm:cofinal}
Under Condition~\ref{cond:geometry}, for every compact $Q\subset M$
there is a connected compact regular domain $\Gamma$ with
$Q\subset\intt\Gamma$ and simply connected exterior $O_\Gamma$.
\end{theorem}

Here a compact regular domain is a compact codimension-zero
submanifold with nonempty smooth boundary.
Theorem~\ref{thm:cofinal} gives simple connectivity at infinity:
loops outside a sufficiently large compact domain can be filled
while avoiding any prescribed compact subset.
Freedman's recognition theorem then yields Theorem~\ref{thm:main}
\cite[Corollary~1.2]{Freedman}.
The argument identifies the homeomorphism type of $M$ and does not
establish a diffeomorphism with standard $\R^4$.

A related result of Chodosh, M\'aximo, and Mukherjee
\cite[Theorem~1.2]{CMM} shows that a compact contractible smooth
4-manifold is homeomorphic to the 4-ball if its interior
admits a complete metric with uniformly positive scalar curvature.
In the present setting, the curvature and volume-growth assumptions
supply control of the end without a compactification hypothesis.
The main construction produces compact cuts whose boundaries admit
positive scalar curvature and which are regular sublevels of global
strictly three-convex Morse exhaustions
(Theorem~\ref{thm:psc-defining}).
This allows boundary geometry and exterior topology to be controlled
for the same domains.

Theorem~\ref{thm:topology} isolates a criterion for such a domain to
have simply connected exterior.
Its assumptions, stated in Condition~\ref{cond:topology}, concern
the boundary, a relative cell model, homology, and the exterior
fundamental group. The group hypothesis only excludes subgroups
isomorphic to $F_2$, and in particular applies to finite or amenable
exterior groups.
The criterion is independent of the original Ricci curvature and
volume-growth hypotheses, and is available for other geometric
constructions whenever all of these cut conditions are verified.
Thus the argument provides both a recognition result and a criterion
that can be applied separately from the analytic construction used
here.

After the preliminaries,
Sections~\ref{sec:analytic}--\ref{sec:finite-exteriors}
construct the analytic data and establish exterior finiteness.
Section~\ref{sec:cuts} constructs the cuts, and
Sections~\ref{sec:topology}--\ref{sec:assembly}
develop the exterior criterion and prove topological recognition.
Detailed proofs appear in the appendices.

\section{Preliminaries}\label{sec:preliminaries}

All manifolds are Hausdorff and second-countable.
Homology and cohomology have integer coefficients unless specified
otherwise. We use $\Delta=\sum_i\partial_i^2$ in Euclidean coordinates
and the curvature convention for which round spheres have positive
sectional curvature.

\subsection{Domains and geometric notation}

\begin{definition}[Compact regular domains]\label{def:ambient}
A compact regular domain $K\subset M$ is a compact embedded
codimension-zero submanifold with nonempty smooth boundary.
\end{definition}

\begin{definition}[The two exteriors]\label{def:exteriors}
Set
\[
 Y_K=\partial K,\qquad
 O_K=M\setminus K,\qquad
 E_K=M\setminus\intt K.
\]
The manifold with boundary $E_K$ has interior $O_K$.
Boundary collars therefore make the inclusion
$O_K\hookrightarrow E_K$ a homotopy equivalence
\cite[Theorems~9.25--9.26]{Lee}.
\end{definition}

\begin{definition}[Differential operators and curvature]
\label{def:operators}
For a Riemannian metric $q$, use
\[
 \Hess_qf=\nabla^qdf,\qquad
 \Delta_qf=\tr_q\Hess_qf,\qquad
 \Scal_q=\tr_q\Ric_q.
\]
For a symmetric bilinear form $S$ and a subspace $P$,
\[
 \tr_P^qS=\sum_{i=1}^{\dim P}S(e_i,e_i),
\]
where $(e_i)$ is a $q$-orthonormal basis of $P$.
Tensor inequalities are inequalities of quadratic forms
\cite[Chapters~2--3]{Petersen}.
\end{definition}

\begin{definition}[Mean curvature and regular levels]
\label{def:mean}
With $\nu_q$ the outward unit normal of $Y_K$, set
\[
 \mathrm{II}_q(V,W)=q(\nabla^q_V\nu_q,W),\qquad
 H_q(Y_K)=\tr_{TY_K}^q\mathrm{II}_q.
\]
If $f$ is smooth near $Y_K$, with $K=\{f\le s\}$ there and
$df\ne0$ on $Y_K$, then $Y_K=\{f=s\}$ locally and
\[
 \nu_q=\frac{\nabla_qf}{|df|_q},\qquad
 H_q(Y_K)=\frac{\tr_{TY_K}^q\Hess_qf}{|df|_q}.
\]
Strict mean convexity means $H_q(Y_K)>0$.
\end{definition}

\begin{definition}[Positive scalar curvature]\label{def:psc}
A manifold admits PSC if it carries a smooth metric with positive
scalar curvature. For $Y_K$, this metric need not be the induced
metric. Uniform positivity means a positive global lower bound.
\end{definition}

\subsection{Conformal formulas}

For $\widetilde q=e^{2h}q$ in dimension four,
\[
 |V|_{\widetilde q}=e^h|V|_q,\qquad
 |\alpha|_{\widetilde q}=e^{-h}|\alpha|_q,
\]
for vectors $V$ and covectors $\alpha$, and
\[
\begin{aligned}
 \Hess_{\widetilde q}f
 &=\Hess_qf-dh\otimes df-df\otimes dh
   +\langle dh,df\rangle_q\,q,\\
 \Ric_{\widetilde q}
 &=\Ric_q-2\Hess_qh+2dh\otimes dh
   -\bigl(\Delta_qh+2|dh|_q^2\bigr)q,\\
 \Scal_{\widetilde q}
 &=e^{-2h}\bigl(\Scal_q-6\Delta_qh-6|dh|_q^2\bigr).
\end{aligned}
\]
For the same outward side of a hypersurface,
\[
 H_{\widetilde q}
 =e^{-h}\bigl(H_q+3\nu_qh\bigr).
\]
The connection, Ricci, and scalar-curvature transformations are
\cite[Exercise~4.7.14(1), (4), and (5)]{Petersen}.
The Hessian formula follows from the connection transformation;
tracing the transformed second fundamental form gives the
mean-curvature formula with our outward-normal convention.

\subsection{Exhaustions and Morse topology}

\begin{definition}[Exhaustions and Morse notation]\label{def:morse}
An exhaustion is a continuous function bounded below with compact
sublevel sets.
For a smooth function $f$, write $\Crit f=\{df=0\}$.
It is Morse when its critical Hessians are nondegenerate;
$\ind_xf$ denotes the number of negative Hessian eigenvalues.
\end{definition}

\begin{definition}[Three-convexity]\label{def:convex}
For a symmetric form $S$ on a 4-dimensional inner product
space $(V,q)$, set
\[
 \lambda_{(3)}^q(S)
 :=\min_{\dim P=3}\tr_P^qS
 =\lambda_1^q(S)+\lambda_2^q(S)+\lambda_3^q(S),
\]
where the eigenvalues are in increasing order.
Strict and weak three-convexity of $f$ mean, respectively,
\[
 \lambda_{(3)}^q(\Hess_qf)>0,\qquad
 \lambda_{(3)}^q(\Hess_qf)\ge0.
\]
\end{definition}

A strictly three-convex Morse function has indices at most two
and strictly mean-convex regular sublevel boundaries.

In a compact 4-dimensional Morse band, an index-$r$ critical
point contributes an $r$-handle and, up to homotopy, a relative
$r$-cell \cite[Theorem~3.2]{Milnor},
\cite[Theorems~3.13--3.14]{MilnorHC}.
Reversing the band changes the handle index to $4-r$, since the
critical Hessian changes sign.
A bound on relative cell dimensions concerns the cells attached
to the boundary, not the dimension of the entire model.

\subsection{Volume comparison and Green potentials}

\begin{definition}[Volume growth]\label{def:AVR}
For a connected complete 4-manifold with $\Ric_g\ge0$,
\[
 \theta:=\AVR(g)
 =\lim_{R\to\infty}
 \frac{\Vol_gB_R^g(p)}{\omega_4R^4},
 \qquad \omega_4=\frac{\pi^2}{2},
\]
is independent of $p$.
Euclidean volume growth means $\theta>0$.
Bishop--Gromov comparison \cite[Lemma~7.1.4]{Petersen} then gives
\[
 \theta\omega_4r^4
 \le\Vol_gB_r^g(x)\le\omega_4r^4
 \qquad(x\in M,\ r>0).
\]
\end{definition}

\begin{definition}[Green potentials]\label{def:greenpotential}
A positive symmetric Green kernel $G_g(x,y)=G_g(y,x)$ for
$-\Delta_g$ is smooth off the diagonal, locally integrable,
and satisfies, for every
$\eta\in C_c^\infty(M)$,
\[
 \int_MG_g(x,y)(-\Delta_g\eta(x))\,d\Vol_g(x)=\eta(y).
\]
For a nonnegative source $f$, write
\[
 P_f(x)=\int_MG_g(x,y)f(y)\,d\Vol_g(y)
\]
where finite; see \cite[Section~5]{LY} and
\cite[Introduction]{LT} for the Green-kernel conventions.
Existence of the positive kernel used here is established in
Proposition~\ref{prop:green}.
\end{definition}

\subsection{Group conventions}

Write $F_r$ for the free group of rank $r$ and $G*_AH$ for
the pushout of specified homomorphisms $A\to G$ and $A\to H$.

A group is residually finite if every nonidentity element
survives under a homomorphism to a finite group.
The regular permutation representation embeds every finite
group into a unitary group.
Consequently, finite-dimensional unitary representations
separate elements of a residually finite group.

\section{Green potential and a three-convex exhaustion}\label{sec:analytic}

We construct a Green potential that yields two complete conformal
metrics: one with uniformly positive scalar curvature and a weighted
Ricci lower bound, and another admitting a strictly three-convex
Morse exhaustion.
The arguments in this section use only completeness, nonnegative
Ricci curvature, and positive asymptotic volume ratio in dimension
four; contractibility is not needed.

\begin{proposition}[Newton potential]\label{prop:green}
Let $(M,g)$ be a connected complete smooth Riemannian 4-manifold
without boundary, with $\Ric_g\ge0$ and $\operatorname{AVR}(g)>0$.
Fix $p\in M$ and put $\rho(x)=1+d_g(p,x)$.
There exist a smooth positive source $f$ and constants $c,C>0$ such that
\[
 \frac{c}{\rho^3}\le f\le\frac{C}{\rho^3}.
\]
The minimal positive Green kernel $G_g$ of $-\Delta_g$ exists, and
the Newton potential
\[
 P_f(x):=\int_M G_g(x,y)f(y)\,d\Vol_g(y)
\]
is finite everywhere and defines a smooth positive function $u=P_f$.
There are constants $L,U,D>0$ such that
\begin{equation}\label{eq:newton}
 -\Delta_g u=f,\qquad
 \frac{L}{\rho}\le u\le\frac{U}{\rho},\qquad
 |du|_g\le\frac{D}{\rho^2}.
\end{equation}
\end{proposition}

These estimates give the following conformal data.

\begin{proposition}[Conformal data]\label{prop:conformal}
Let $u$ be the potential constructed in Proposition~\ref{prop:green}.
There are constants $m,B>0$ such that, with
\[
 \varphi=\log u,\qquad \gbar=u^2g=e^{2\varphi}g,
\]
the metric $\gbar$ is complete, $\sup_M\varphi<\infty$, and
$-\varphi$ is an exhaustion. Moreover,
\begin{equation}\label{eq:weighted}
\begin{gathered}
 -\Delta_g u\ge m u^3,\qquad
 |d\varphi|_{\gbar}^2\le B,\qquad
 \Scal_{\gbar}\ge6m,\\
 \Ric_{\gbar}
 +2\Hess_{\gbar}\varphi
 +d\varphi\otimes d\varphi
 \ge m\gbar.
\end{gathered}
\end{equation}
\end{proposition}

The same potential also yields a strictly three-convex Morse
exhaustion for a second complete conformal metric.
The resulting index bound will control the dimensions of the
relative cells in the exterior of a cut.

\begin{theorem}[Three-convex Morse exhaustion]\label{thm:joint}
For the potential $u$ in Proposition~\ref{prop:conformal}, set
\[
 \psi=\frac{\varphi}{3}=\frac{\log u}{3},
 \qquad q=u^{2/3}g=e^{2\psi}g.
\]
Then $(M,q)$ is complete and admits a smooth Morse exhaustion $F$
satisfying
\[
 \lambda_{(3)}^q(\Hess_qF)>0
 \qquad\text{on }M.
\]
In particular,
\[
 \ind_xF\le2\qquad(x\in\Crit F).
\]
\end{theorem}

The choice of exponent in $q=u^{2/3}g=e^{-4\varphi/3}\gbar$
connects the weighted Ricci estimate to three-convexity.
To see this, let $b$ be a ray Busemann function for $\gbar$,
and let $h$ be one of the smooth local lower supports constructed
in Lemma~\ref{lem:weighted-busemann}.
At the contact point, $|dh|_{\gbar}=1$ and
$\Hess_{\gbar}h(\nabla_{\gbar}h,\nabla_{\gbar}h)=0$.
The conformal Hessian formula therefore gives
\[
 \tr_{\ker dh}^{q}\Hess_qh
 =e^{4\varphi/3}
  \left(\Delta_{\gbar}h
       -2\langle d\varphi,dh\rangle_{\gbar}\right).
\]
An index-form estimate using \eqref{eq:weighted} makes the
expression in parentheses positive. Thus the exponent matches
the drift coefficient supplied by the weighted Ricci bound.

Positivity on $\ker dh$ alone does not yet control every
three-plane. An increasing convex reparametrization $\chi$
contributes the positive term $\chi''(h)\,dh\otimes dh$ to the Hessian;
choosing $\chi''/\chi'$ sufficiently large on successive compact
sets gives positivity on all three-planes.
The Busemann function is proper by the annular diameter estimate
in Lemma~\ref{lem:proper-busemann}. Shen's approximation lemma
then produces the smooth strictly three-convex Morse exhaustion
\cite[Section~3.1, Lemma~5]{Shen}.
The support construction and the locally uniform estimates needed
for this approximation are proved in Appendix~\ref{app:morse}.

\begin{remark}
For the Euclidean model on $\R^4\setminus\{0\}$, take
$u=r^{-1}$, where $r=|x|$.
With $t=\log r$ and $s=\frac32r^{2/3}$, the conformal metrics become
\[
 \gbar=dt^2+g_{S^3},\qquad \Scal_{\gbar}=6,\qquad
 q=ds^2+\left(\frac{2s}{3}\right)^2g_{S^3}.
\]
Moreover,
\[
 \Hess_q t=\frac{3}{2s^2}(q-2\,ds\otimes ds),\qquad
 \lambda_{(3)}^q(\Hess_qt)=\frac{3}{2s^2}>0.
\]
This model illustrates the different roles of the two conformal
metrics and is considered only away from the origin.
\end{remark}

\section{Finite exterior fundamental groups}\label{sec:finite-exteriors}

We retain the complete metric $q$ of Theorem~\ref{thm:joint}.
Every strictly three-convex Morse exhaustion for $q$ yields a
cofinal family of connected mean-convex cores with uniformly
finite exterior fundamental groups.
The proof combines the annular estimate of
\cite[Theorem~4.8]{CMS} with the exterior surjections supplied
by the Morse index bound; see Appendix~\ref{app:finite-exteriors}.

We state the result for an arbitrary exhaustion $H$, since
Section~\ref{sec:cuts} uses both the exhaustion $F$ constructed
above and the global defining functions $G$ associated with
the chosen cuts.

\begin{theorem}[Finite exterior fundamental groups]\label{thm:cores}
Under Condition~\ref{cond:geometry}, let $q$ be the complete metric
fixed above, and let $H$ be any smooth strictly three-convex Morse
exhaustion for $q$.
There are strictly increasing regular values $t_i\to\infty$,
with $H(p)<t_1$. Let $D_i$ be the component of $\{H\le t_i\}$
containing $p$. Each $D_i$ is a connected compact regular domain,
and both $Y_{D_i}$ and $O_{D_i}$ are path connected.
There is a single constant $N_H\in\N$ such that
\begin{equation}\label{eq:cores}
\begin{aligned}
 &D_i\subset\intt D_{i+1}\quad(i\ge1),\qquad
   \bigcup_{i\ge1}\intt D_i=M,\\
 &H_q(Y_{D_i})>0\quad(i\ge1),\\
 &|\pi_1(O_{D_i})|\le N_H\quad(i\ge1),\\
 &\pi_1(O_{D_j})\twoheadrightarrow\pi_1(O_{D_i})
   \quad(j\ge i\ge1).
\end{aligned}
\end{equation}
The bound $N_H$ can be chosen to depend only on $(M,g)$,
independently of $H$.
The displayed maps are induced by the actual inclusions
and compatible basepoint paths.
\end{theorem}

For the exhaustion $F$ of Theorem~\ref{thm:joint}, we write
$s_i$ and $C_i$ for these levels and cores; they are fixed
from now on.

\section{PSC cuts}\label{sec:cuts}

We construct compact domains whose boundaries admit PSC and which
are regular sublevels of global strictly three-convex Morse
exhaustions. The identity $\Gamma=\{G\le0\}$ allows the exterior
finiteness theorem to be applied to larger cores of the same
exhaustion.

Throughout this section, fix the conformal and Morse data and the
compact regular domains $C_i$ from
Theorems~\ref{thm:joint}--\ref{thm:cores}.
In particular, the metric $q$ remains fixed.

\subsection{Boundary connectedness and fundamental groups}

We first record the boundary theorem used to obtain connectedness
and surjectivity on fundamental groups.

\begin{theorem}[Spectral boundary theorem]\label{input:spectral}
Let $(M,g)$ be a smooth Riemannian 4-manifold, and let
$K\subset M$ be a connected compact regular domain.
Suppose $v>0$ is smooth near $K$ and
\[
\begin{aligned}
 (-\Delta_gv)g+v\Ric_g&>0 &&\text{on }K,\\
 H_g(Y_K)+v^{-1}\nu_gv&\ge0 &&\text{on }Y_K.
\end{aligned}
\]
Then $Y_K$ is connected, and inclusion induces a surjection
\[
 \pi_1(Y_K)\twoheadrightarrow\pi_1(K).
\]
This is the compact case of \cite[Theorem~2]{ALS}, with $n=4$
and $\gamma=1$. The tensor inequality is equivalent to
$-\Delta_gv+v\lambda_{\min}^g(\Ric_g)>0$; compactness of $K$
supplies completeness and a compact boundary component, and
$1\le(4-1)/(4-2)=3/2$ is in the permitted range.
\end{theorem}

\begin{corollary}[Mean-convex boundary groups]\label{cor:boundaryepi}
For the fixed conformal data, every connected compact regular
domain $K$ with $H_q(Y_K)>0$ has connected boundary, and inclusion
induces a surjection
\[
 \pi_1(Y_K)\twoheadrightarrow\pi_1(K).
\]
\end{corollary}

\begin{proof}
We have
\[
 (-\Delta_gu)g+u\Ric_g\ge mu^3g>0.
\]
Since $q=u^{2/3}g$, the conformal mean-curvature formula gives
\begin{equation}\label{eq:meanconformal}
 H_q(Y_K)
 =u^{-1/3}\bigl(H_g(Y_K)+\nu_g\log u\bigr).
\end{equation}
Apply Theorem~\ref{input:spectral} with $v=u$.
\end{proof}

\subsection{Forcing and an outermost PSC boundary}

The forcing function provides the boundary barriers and the
positivity needed to deduce PSC from stability.

\begin{proposition}[Late positive forcing]\label{prop:forcing}
For every $i$, there are $j>i$ and a positive function
$h\in C^\infty(M)$ such that, with
\[
 A=C_i,\qquad B=C_j,\qquad k=u^{2/3}h,
\]
we have
\begin{equation}\label{eq:forcing}
\begin{gathered}
 A\subset\intt B,\qquad
 h^2-2|dh|_{\gbar}\ge-m,\\
 H_q(Y_A)<k,\qquad H_q(Y_B)>k.
\end{gathered}
\end{equation}
Moreover,
\begin{equation}\label{eq:DEC}
 \Scal_q+\frac43k^2-2|dk|_q>0.
\end{equation}
\end{proposition}

\begin{definition}[Weakly outermost compact regular domain]
\label{def:outermost}
Let $(M,q)$ be a smooth Riemannian 4-manifold, let $B$ be a
compact regular domain, and let $k\in C^\infty(M)$ be positive on $B$.
A compact regular domain $T\subset\intt B$ is
\emph{weakly outermost for $H_q\le k$ in $B$} if
$H_q(Y_T)\le k$ and the following property holds:
every compact regular domain $L$ satisfying
\[
 T\subseteq L\subset\intt B,\qquad H_q(Y_L)\le k
\]
must equal $T$.
\end{definition}

We use the following existence theorem both for the PSC boundary
and for the hull capture argument below.

\begin{theorem}[Outermost prescribed-mean-curvature domain]
\label{input:eichmair}
Let $(M,q)$ be a complete connected oriented Riemannian
4-manifold, and let $A\subset\intt B$ be compact regular domains.
Let $\psi\in C^\infty(M)$ be positive on $B$, and suppose
\[
 H_q(Y_A)<\psi,\qquad H_q(Y_B)>\psi.
\]
Then there is a compact regular domain $E$ with
\[
 A\subset\intt E,\qquad E\subset\intt B,\qquad H_q(Y_E)=\psi,
\]
which contains every compact regular domain $D$ satisfying
\[
 A\subseteq D\subset\intt B,\qquad H_q(Y_D)\le\psi.
\]
This follows from \cite[Theorem~5.1]{EichmairHorizons} with
initial-data tensor $(\psi/3)q$; see Appendix~\ref{app:union}.
\end{theorem}

\begin{proposition}[Outermost PSC boundary]\label{prop:outerpsc}
For $A,B,k$ in Proposition~\ref{prop:forcing}, there is a compact
regular domain $T$, weakly outermost for $H_q\le k$ in $B$, such that
\[
 A\subset\intt T,\qquad T\subset\intt B,\qquad H_q(Y_T)=k.
\]
Its boundary $Y_T$ admits PSC.
\end{proposition}

\begin{corollary}[The component containing the inner core]
\label{cor:psccomponent}
Let $K$ be the connected component of $T$ containing $A$.
Then $K$ is a compact regular domain with
\[
 A\subset\intt K,\qquad H_q(Y_K)>0.
\]
The boundary $Y_K$ is connected and admits PSC, and inclusion induces
\[
 \pi_1(Y_K)\twoheadrightarrow\pi_1(K).
\]
\end{corollary}

\begin{proof}
The boundary $Y_K$ is a union of components of $Y_T$, and
$H_q(Y_K)=k>0$. Apply Corollary~\ref{cor:boundaryepi}.
\end{proof}

\subsection{Inward collars}

A sufficiently small inward displacement preserves strict mean
convexity and the diffeomorphism type of the boundary.
We use this freedom in the hull capture argument.

\begin{lemma}[Inward collar cut]\label{lem:psc-collar}
Let $K$ be the component in Corollary~\ref{cor:psccomponent}.
There is a connected compact regular domain $\Gamma$ such that
\[
 A\subset\intt\Gamma,\qquad
 \Gamma\subset\intt K,\qquad H_q(Y_\Gamma)>0,
\]
and $Y_\Gamma$ is diffeomorphic to $Y_K$.
Consequently, $Y_\Gamma$ is connected and admits PSC, and inclusion induces
\[
 \pi_1(Y_\Gamma)\twoheadrightarrow\pi_1(\Gamma).
\]
\end{lemma}

\subsection{A global defining Morse level for the cut}

For a nonempty compact set $K$, define its weak three-convex hull by
\[
 \widehat K^{\,q}_3
 =\left\{x\in M:
 \begin{array}{l}
  v(x)\le\max_Kv\text{ for every }v\in C^\infty(M)\\
  \text{with }\lambda_{(3)}^q(\Hess_qv)\ge0
 \end{array}
 \right\}.
\]
This is the ambient smooth $p$-convex hull of
\cite[Definition~3.1]{HarveyLawson}, with $p=3$ and the Riemannian
Hessian convention of \cite[Definition~2.2$'$]{HarveyLawson}.
Hull equality means that every point outside $K$ can be separated
from $K$ by a smooth weakly three-convex function.
It is an ambient separation property; three-convexity of the
interior alone does not assert this equality.

The following maximum principle gives a useful way to detect
this separation property without constructing the separators
individually.

\begin{definition}[Relative maximum principle]\label{def:bypass-hull}
For a closed set $S$ and an open set $\Omega$, write
$\mathsf{MP}_q(S;\Omega)$ if every smooth local test function $\psi$
having a local maximum on $S$ at $x\in S\cap\Omega$ admits
some three-plane $P\subset T_xM$ such that
\[
 \tr_P^q\Hess_q\psi(x)\le0.
\]
Local tests may equivalently be extended to global smooth functions
by a cutoff that is identically one near the contact point.
\end{definition}

\begin{lemma}[Direct compact capture]\label{lem:bypass-capture}
Let $A,B,k,T$ be the barriers and outermost domain of
Proposition~\ref{prop:outerpsc}, and let $K$ be its component
containing $A$.
Let $\Gamma\subset\intt K$ be an inward collar cut.
Suppose $M$ has a smooth strictly three-convex exhaustion $F$,
and choose a smooth collar defining function $r$ and a number
$c<0$ satisfying
\[
 \Gamma=\{x\in K:r(x)\le c\},
\]
and, throughout $K\setminus\Gamma$,
\[
 dr\ne0,\qquad
 \tr_{(\nabla_qr)^\perp}^q\Hess_qr>0.
\]
Every compact $S\subset B$ satisfying
$\mathsf{MP}_q(S;M\setminus\Gamma)$ is contained in $\Gamma$.
\end{lemma}

The key point is to fix a protection radius using only the
positive lower bound for $k$, the ambient geometry on $B$, and
the separation of $\Gamma$ from $Y_T$; see
Lemma~\ref{lem:bypass-protection}.
One then enlarges $T$ to contain $S$ and increases the forcing
only in a smaller neighborhood of $S$.
The prescribed-mean-curvature theorem gives an enclosing domain
$E$ whose boundary stays outside that neighborhood by the same
protection estimate.
Thus $H_q(Y_E)=k$ on the new boundary, and the original
outermostness forces $E=T$.
The exhaustion excludes other components of $T$, and a convex
reparametrization of the collar coordinate excludes
$S\setminus\Gamma$.
The full proof, including the uniform protection estimate, is
given in Appendix~\ref{app:union}.

The weak hull satisfies the relative maximum principle, and
the exhaustion makes it compact.
For the cuts constructed here, an enclosing regular level of
$F$ also places the hull inside $B$.
Lemma~\ref{lem:bypass-capture} therefore gives the following
ambient hull equality.

\begin{lemma}[Hull capture for inward cuts]\label{lem:bypass-hull-body}
The inward collar construction in Lemma~\ref{lem:psc-collar}
can be made so that
\[
 \widehat\Gamma^{\,q}_3=\Gamma.
\]
\end{lemma}

Together with strict mean convexity, this separation property
allows the collar defining function to extend to a global
strictly three-convex exhaustion: convex reparametrization makes
the local defining function strict, while hull separators and
smooth maxima join it to an exhaustion away from the boundary.
A fine Morse approximation fixed near the zero level then
preserves the chosen domain.

\Needspace{20\baselineskip}
\begin{theorem}[PSC cuts with global defining Morse levels]
\label{thm:psc-defining}
Under Condition~\ref{cond:geometry}, retain the fixed data of
Theorems~\ref{thm:joint}--\ref{thm:cores}.
For every compact $Q\subset M$, there are a connected compact
regular domain $\Gamma$ and a smooth strictly three-convex Morse
exhaustion $G$ for $q$ with
\[
 Q\subset\intt\Gamma,\qquad H_q(Y_\Gamma)>0.
\]
The boundary $Y_\Gamma$ is connected and admits PSC, and inclusion induces
\[
 \pi_1(Y_\Gamma)\twoheadrightarrow\pi_1(\Gamma).
\]
The domain and its boundary are given by
\[
 \Gamma=\{G\le0\},\qquad Y_\Gamma=\{G=0\},
\]
with $dG\ne0$ on $Y_\Gamma$. In addition,
\[
\begin{aligned}
 \lambda_{(3)}^q(\Hess_qG)&>0 &&\text{on }M,\\
 \ind_xG&\le2 &&\text{for }x\in\Crit G.
\end{aligned}
\]
The domain $\Gamma$ and the function $G$ may depend on $Q$,
while $q$ remains fixed.
\end{theorem}

The proofs of the hull capture lemma and the global defining-level
theorem are given in Appendix~\ref{app:union}.

\section{From a Morse band to a simply connected exterior}
\label{sec:topology}

This section gives a topological criterion for a PSC cut to have
simply connected exterior. A Morse band supplies a relative cell
model and, in the geometric application, a finite exterior group.
The criterion itself only requires that the exterior group contain
no subgroup isomorphic to the free group $F_2$.
The PSC boundary decomposition, a cellular form of the
Gerstenhaber--Rothaus theorem, and a smooth filling obstruction
then determine the boundary-group map needed for van Kampen.

The criterion is stated in Theorem~\ref{thm:topology} under
Condition~\ref{cond:topology}, independently of the geometric
construction in Section~\ref{sec:cuts}.
Detailed proofs are collected in Appendix~\ref{app:topology}.

For a group $G$, write $H_i(G;\Z)=H_i(BG;\Z)$, where $BG$ is
a connected CW complex with fundamental group $G$ and contractible
universal cover. In particular, $H_1(G;\Z)=G/[G,G]$.
For a space $X$, the groups $H_2(X)$ and $H_2(\pi_1(X);\Z)$
need not agree; their comparison below uses the Hopf homology
surjection \cite[Section~0]{Brown}.

\subsection{Topological data supplied by a Morse band}

\begin{definition}[Mean-convex band Morse data]
\label{def:band}
Let $(M,q)$ be a smooth Riemannian 4-manifold, and let
$K\subset\intt B$ be compact regular domains.
Write $\mathsf{Band}_q(K,B)$ if there is $a\in C^\infty(M)$
whose critical points in $B\setminus\intt K$ have nondegenerate
Hessians. Its boundary values and, for every $x$ in this band,
its pointwise conditions must satisfy
\begin{equation}\label{eq:band}
\begin{aligned}
 a|_{Y_K}&=0,\qquad a|_{Y_B}=1,\\
 \ind_xa&\le2 &&\text{if }da_x=0,\\
 \tr_{(\nabla_q a)^\perp}^{q}\Hess_q a(x)&>0
     &&\text{if }da_x\ne0.
\end{aligned}
\end{equation}
No range restriction is imposed on the interior values of $a$.
\end{definition}

\begin{lemma}[Boundary orientation of band data]
\label{lem:bandorientation}
Suppose $\mathsf{Band}_q(K,B)$ and
$H_q(Y_K),H_q(Y_B)>0$.
Every function $a$ realizing these band data satisfies
\[
 da\ne0\quad\text{on }Y_K\cup Y_B,
\]
and $\nabla_q a$ points outward from both $K$ and $B$.
\end{lemma}

\begin{definition}[Relative cell model]
\label{def:CW}
Fix a finite CW structure on $Y_K$.
A \emph{relative $0,1,2$-cell model} of $(E_K,Y_K)$ consists of
a CW complex $X$, a cell-preserving identification $j$ of $Y_K$
with a subcomplex of $X$, and a homotopy equivalence
$e:X\to E_K$ such that
\[
 e\circ j=(Y_K\hookrightarrow E_K).
\]
The complex $X$ is obtained from $j(Y_K)$ by attaching relative
cells of dimensions $0,1,2$ and need not be finite.
\end{definition}

\begin{proposition}[Relative topology of a Morse cut]
\label{prop:transfer}
Under Condition~\ref{cond:geometry}, let $H$ and $D_i$ be as in
Theorem~\ref{thm:cores}.
Suppose $K\subset\intt D_j$ is a connected compact regular domain
with connected boundary, $H_q(Y_K)>0$, and
$\mathsf{Band}_q(K,D_j)$.
Then $O_K$ is path connected, inclusion induces a surjection
\[
 \pi_1(O_{D_j})\twoheadrightarrow\pi_1(O_K),
\]
and
\begin{equation}\label{eq:transfer}
 H_\ell(M,K)=0\quad(\ell\ge3),\qquad
 H_2(M,K)\text{ is torsion-free}.
\end{equation}
Every finite CW structure on $Y_K$ extends to a relative cell
model as in Definition~\ref{def:CW}.
\end{proposition}

\begin{corollary}[Homology and exterior finiteness]
\label{cor:transferdata}
Under the hypotheses of Proposition~\ref{prop:transfer},
$H_1(K)$ is torsion-free and
\[
 |\pi_1(O_K)|\le N_H,\qquad H_2(K)=H_3(K)=0,
\]
where $N_H$ is the bound in Theorem~\ref{thm:cores}.
\end{corollary}

\begin{proof}
The exterior group is a quotient of $\pi_1(O_{D_j})$.
Contractibility of $M$ and the exact sequence of $(M,K)$ give
\[
 \begin{gathered}
  H_2(K)\cong H_3(M,K)=0,\\
  H_3(K)\cong H_4(M,K)=0,\\
  H_1(K)\cong H_2(M,K).
 \end{gathered}
\]
The assertions follow from Proposition~\ref{prop:transfer}.
\end{proof}

We collect these consequences with the boundary data needed
for the topological argument.

\begin{condition}[Topological cut hypotheses]
\label{cond:topology}
Let $M$ be a contractible smooth 4-manifold, and let
$K\subset M$ be a connected compact regular domain with connected
boundary.
Write $Y=Y_K$, $O=O_K$, and $E=E_K$.
Assume that $Y$ admits PSC, $O$ is path connected, and
$\pi_1(O)$ contains no subgroup isomorphic to $F_2$.
Assume also that $H_1(K)$ is torsion-free and
\begin{equation}\label{eq:topdata}
 \pi_1(Y)\twoheadrightarrow\pi_1(K),\qquad H_2(K)=H_3(K)=0.
\end{equation}
The displayed boundary-group map is induced by inclusion.
Assume also that the simultaneous boundary data of
Proposition~\ref{input:normalform} can be chosen so that the resulting
finite CW structure on $Y$ extends to a relative cell model as in
Definition~\ref{def:CW}. Only this chosen CW structure is required.
\end{condition}

The exclusion of $F_2$ holds, in particular, for finite groups and
for amenable groups. In the geometric application,
Corollary~\ref{cor:transferdata} supplies finiteness, while
Proposition~\ref{prop:transfer} supplies a relative model for every
finite boundary CW structure and hence for the one required here.

\begin{lemma}[Boundary homology maps]
\label{lem:boundaryH}
Let $M$ be a contractible smooth 4-manifold, and let
$K\subset M$ be a connected compact regular domain with connected
boundary.
If $H_2(K)=H_3(K)=0$ and $H_1(K)$ is torsion-free, inclusion induces
\begin{equation}\label{eq:boundaryH}
\begin{gathered}
 H_1(Y_K)\xrightarrow{\cong}H_1(K),\qquad
 H_2(Y_K)\xrightarrow{\cong}H_2(E_K),\\
 H_1(E_K)=0.
\end{gathered}
\end{equation}
In particular, $H_1(Y_K)$ is torsion-free.
\end{lemma}

\subsection{PSC boundary cells and the auxiliary filling}

The next construction records two consequences of the same PSC
connected-sum decomposition. An acyclic subcomplex retains the
groups of the homology-sphere summands, while boundary surgeries
remove the $S^1\times S^2$ summands.
The specified connecting homomorphism in \eqref{eq:actualdelta}
relates the homology of the resulting filling to that of $K$.

\begin{proposition}[PSC boundary cells and auxiliary filling]
\label{input:normalform}
Let $K$ be a connected compact regular domain in a simply
connected smooth 4-manifold.
Suppose $Y=\partial K$ is connected, admits PSC, and
$H_1(Y)$ is torsion-free.
Let $P$ denote the Poincar\'e homology sphere.

The following data can be constructed simultaneously.

\begin{enumerate}
\item
A finite CW structure on $Y$ with one three-cell and no higher
cells, together with a finite connected acyclic subcomplex
\[
 A\subset Y^{(2)},\qquad \dim A\le2.
\]
There is a selected family of open two-cells
$e^2_1,\ldots,e^2_s$ whose attaching maps are constant at points
of $A$. Every other two-cell belongs to $A$.
The group $\pi_1(A)$ is residually finite.

\item
A connected compact regular domain $D$ in an auxiliary smooth
4-manifold and a continuous embedding $f:K\hookrightarrow D$
such that, for some integers $r,v\ge0$,
\[
 \begin{gathered}
  \partial D\cong\mathop{\#}_{\ell=1}^{v}(\pm P),\qquad
  \pi_1(\partial D)\cong\pi_1(A),\\
  \pi_1(Y)\cong\pi_1(\partial D)*F_r,\qquad
  H_1(Y)\cong\Z^r.
 \end{gathered}
\]
Each summand carries one of the two orientations of $P$,
and the empty sum is $S^3$.

\item
Identifying $K$ with $f(K)$ when writing the pair $(D,K)$,
\[
 H_1(D,K)=0,\qquad H_2(D,K)\cong\Z^r.
\]
For the connecting homomorphism
$\delta:H_2(D,K)\to H_1(K)$, there are bases
$(e_1,\ldots,e_r)$ of $H_2(D,K)$ and
$(b_1,\ldots,b_r)$ of $H_1(Y)$ such that
\begin{equation}\label{eq:actualdelta}
 \delta e_i=(Y\hookrightarrow K)_*b_i\qquad(1\le i\le r).
\end{equation}
\end{enumerate}
\end{proposition}

The construction uses the PSC classification
\cite[Proposition~2.2 and Corollary~2.3]{CMM},
the spine existence and regular-neighborhood results
\cite[Theorems~2.5 and~2.1]{Spines},
and the boundary surgery modeled by $W_2$ in
\cite[Example~8.26(i)]{Ranicki}, with $m=3$ and $i=1$.

\subsection{Injectivity of the acyclic boundary subcomplex}

The following theorem is a cellular reformulation of the
Gerstenhaber--Rothaus theorem \cite[Theorem~1(ii)]{GR}.
The corresponding solvability statement under a vanishing
second-homology hypothesis, with solutions in the same unitary
group, is also stated in \cite[Lemma~2.1(1)]{NitscheThom}.
Residual finiteness converts this extension property into
injectivity of the subcomplex group. The exterior criterion
combines this extension property with the boundary cells and
the smooth filling obstruction.

\begin{theorem}[A cellular form of the Gerstenhaber--Rothaus theorem]
\label{thm:unitary}
Let $(Z,A)$ be a finite CW pair, with $A$ and $Z$ connected.
Suppose all relative cells have dimensions $0,1,2$ and
$H_2(Z,A)=0$.
For every integer $d\ge1$ and every representation
$\rho:\pi_1(A)\to U(d)$, there is a representation
$\widetilde\rho:\pi_1(Z)\to U(d)$ satisfying
\[
 \widetilde\rho\circ i_*=\rho,\qquad i:A\hookrightarrow Z.
\]
If $\pi_1(A)$ is residually finite, then $i_*$ is injective.
\end{theorem}

\begin{corollary}[An acyclic subcomplex of a two-complex]
\label{cor:spine}
Let $Z$ be a connected, possibly infinite CW complex of dimension
at most two, with $H_2(Z)=0$.
Let $A\subset Z$ be a finite connected subcomplex such that
$H_1(A)=0$ and $\pi_1(A)$ is residually finite.
Then inclusion induces an injection
\[
 \pi_1(A)\hookrightarrow\pi_1(Z).
\]
In particular, if $\pi_1(Z)$ is finite, then $\pi_1(A)$ is finite.
\end{corollary}

The selected sphere cells in Proposition~\ref{input:normalform}
account for the second homology of the exterior.
Removing them produces the two-complex to which
Corollary~\ref{cor:spine} applies.

\begin{lemma}[The deleted two-complex]
\label{lem:deleted}
Under Condition~\ref{cond:topology}, choose the boundary CW
structure and the simultaneous data of
Proposition~\ref{input:normalform} for which the required relative
cell model exists. Let $A$ be its acyclic subcomplex.
Choose such a relative cell model $(X,j,e)$,
and set
\[
 Z:=X^{(2)}\setminus
       \bigcup_{\alpha=1}^{s}j(e^2_\alpha).
\]
Then $Z$ is a connected CW complex of dimension at most two
containing $j(A)$ as a subcomplex, and
\[
 H_2(Z)=0,\qquad
 \pi_1(Z)\cong\pi_1(E_K).
\]
Consequently, inclusion and the comparison maps induce an injection
\[
 \pi_1(A)\hookrightarrow\pi_1(E_K)\cong\pi_1(O_K).
\]
\end{lemma}

\subsection{From the auxiliary filling to a free boundary group}

The remaining obstruction comes from smooth 4-dimensional
topology. For an oriented smooth manifold, the obstruction to
a spin structure is $w_2(TD)\in H^2(D;\Z/2)$
\cite[Sections~2.4 and~3.2]{FK}.
We write $\sigma(D)$ for the signature of its intersection form;
in particular, $H_2(D;\R)=0$ implies $\sigma(D)=0$
\cite[Section~6.4]{Ranicki}.

\begin{theorem}[Rokhlin]
\label{input:rokhlin}
Let $D$ be a compact smooth spin 4-manifold whose boundary is
an integral homology sphere $\Sigma$.
Then
\[
 \sigma(D)\equiv8\mu(\Sigma)\pmod{16},
\]
where $\mu(\Sigma)\in\Z/2$ is the Rokhlin invariant.
Moreover, $\mu(P)=\mu(-P)=1$
\cite[Sections~3.8--3.9]{FK}.
\end{theorem}

\begin{proposition}[The free boundary group]
\label{prop:freeboundary}
Under Condition~\ref{cond:topology}, choose the simultaneous
data $A,D,r,v$ of Proposition~\ref{input:normalform} compatible
with the relative cell model in that condition.
Then $v=0$, $\partial D\cong S^3$, and
\begin{equation}\label{eq:spinefilling}
\begin{gathered}
 \pi_1(A)=1,\qquad \pi_1(\partial D)=1,\qquad \pi_1(Y)\cong F_r,\\
 H_1(D)=H_2(D)=0.
\end{gathered}
\end{equation}
\end{proposition}

Indeed, $\pi_1(A)$ is a free product of $v$ copies of the binary
icosahedral group $\pi_1(P)$ and injects into $\pi_1(O_K)$.
Two or more such factors contain $F_2$, so the group hypothesis
forces $v\le1$. The auxiliary filling has vanishing first and
second homology, and Rokhlin's theorem excludes the remaining
case $\partial D=\pm P$. This is the only step at which the
restriction on the exterior group is used.

\subsection{From the boundary group to the exterior}

Once the boundary group is free, its inclusion-induced epimorphism
onto $\pi_1(K)$ is determined by the homology already established.

\begin{lemma}[A free epimorphism with trivial homological kernel]
\label{lem:freeepi}
Let $r\ge0$, and let $a:F_r\twoheadrightarrow G$ be a surjective
group homomorphism.
If $H_1(a;\Z)$ is injective and $H_2(G;\Z)=0$,
then $a$ is an isomorphism.
\end{lemma}

\begin{theorem}[A criterion for simply connected exteriors]
\label{thm:topology}
Under Condition~\ref{cond:topology}, $O_K$ is simply connected.
\end{theorem}

\begin{proof}
By Proposition~\ref{prop:freeboundary}, $\pi_1(Y)\cong F_r$.
For the inclusion-induced epimorphism
\[
 a:\pi_1(Y)\twoheadrightarrow\pi_1(K),
\]
Lemma~\ref{lem:boundaryH} makes $H_1(a;\Z)$ an isomorphism.
The Hopf homology surjection \cite[Section~0]{Brown} gives
\[
 0=H_2(K)\twoheadrightarrow H_2(\pi_1(K);\Z).
\]
Lemma~\ref{lem:freeepi} therefore shows that $a$ is an isomorphism.

The collared van Kampen pushout now gives
\begin{equation}\label{eq:endgame}
 1=\pi_1(M)
 \cong\pi_1(K)*_{\pi_1(Y)}\pi_1(E_K)
 \cong\pi_1(E_K).
\end{equation}
Since $O_K\simeq E_K$ and $O_K$ is path connected,
$O_K$ is simply connected.
\end{proof}

\section{Cofinality and recognition}\label{sec:assembly}

We now apply the topological criterion to the cuts from
Section~\ref{sec:cuts}. The defining exhaustion of each cut
supplies the Morse band and finite exterior group needed for
this application.

\begin{proposition}[The defining level supplies a Morse band]
\label{prop:global-band}
Let $\Gamma,G$ be as in Theorem~\ref{thm:psc-defining}, and apply
Theorem~\ref{thm:cores} to $H=G$.
For a sufficiently large positive level $s$ in the resulting
regular sequence, let $B$ be the component of $\{G\le s\}$
containing $p$.
Then
\[
 \Gamma\subset\intt B,\qquad H_q(Y_B)>0,\qquad
 |\pi_1(O_B)|<\infty.
\]
Moreover, $\mathsf{Band}_q(\Gamma,B)$ holds and is realized by
\[
 a=\frac{G}{s},\qquad
 a|_{Y_\Gamma}=0,\qquad a|_{Y_B}=1.
\]
Consequently, $\Gamma$ satisfies Condition~\ref{cond:topology}.
\end{proposition}

The proof is given in Appendix~\ref{app:assembly}.

\begin{proof}[Proof of Theorem~\ref{thm:cofinal}]
Given a compact set $Q\subset M$, Theorem~\ref{thm:psc-defining}
gives $\Gamma,G$ with $Q\subset\intt\Gamma$.
Proposition~\ref{prop:global-band} verifies
Condition~\ref{cond:topology} for this $\Gamma$.
Theorem~\ref{thm:topology} therefore makes $O_\Gamma$
simply connected.
\end{proof}

\begin{definition}[Simple connectivity at infinity]
\label{def:SCI}
A noncompact manifold $M$ is \emph{simply connected at infinity}
if, for every compact $Q\subset M$, there is a compact
$L\supset Q$ such that every loop in $M\setminus L$ extends
to a disk in $M\setminus Q$.
\end{definition}

\begin{proposition}\label{prop:SCI}
Under Condition~\ref{cond:geometry}, $M$ is simply connected
at infinity.
\end{proposition}

\begin{proof}
Given compact $Q$, choose $\Gamma$ from
Theorem~\ref{thm:cofinal} and take $L=\Gamma$.
Every loop in $M\setminus L=O_\Gamma$ contracts there,
hence in $M\setminus Q$.
These connected complements form a cofinal family of
neighborhoods of infinity, so $M$ also has one end.
\end{proof}

\begin{theorem}[Contractible open 4-manifold recognition]
\label{thm:Freedman}
A contractible open topological 4-manifold that is simply
connected at infinity is homeomorphic to $\R^4$
\cite[Corollary~1.2, p.~366]{Freedman}.
\end{theorem}

\begin{proof}[Proof of Theorem~\ref{thm:main}]
Under Condition~\ref{cond:geometry}, $M$ is an open contractible
4-manifold; in particular, $\pi_1(M)=1$ and $H_2(M;\Z)=0$.
Proposition~\ref{prop:SCI} and Theorem~\ref{thm:Freedman}
therefore give $M\cong_{\mathrm{top}}\R^4$.
\end{proof}

\section{Conclusion}\label{sec:conclusion}

The argument separates the geometric construction of cuts from
the topological recognition of their exteriors.
The analytic construction in Section~\ref{sec:analytic} does not
require contractibility. In the contractible setting,
Lemma~\ref{lem:bypass-capture} gives the compact capture needed
to identify an ambient three-convex hull, and
Theorem~\ref{thm:psc-defining} realizes the resulting PSC cut as
a regular sublevel of a global strictly three-convex Morse exhaustion.
This realization makes the boundary geometry and the exterior
Morse topology available for the same domain.

Theorem~\ref{thm:topology} can be used independently whenever its
boundary, homological, and relative cell hypotheses hold and
the exterior group contains no subgroup isomorphic to $F_2$.
It therefore includes finite and amenable exterior groups.
The cellular Gerstenhaber--Rothaus theorem supplies the algebraic
input; its combination with the boundary cell structure and the
smooth filling obstruction is what yields the exterior criterion.
The recognition obtained here is topological. A diffeomorphism
conclusion would require additional control of the smooth structure.

\section*{AI Usage}

The authors used ChatGPT (OpenAI; primarily GPT-6.1 Sol, with occasional use of GPT-6.1 Terra) and Claude (Anthropic; Opus 5.5).
These tools assisted with literature searches, mathematical exploration, checking arguments and calculations, and refining constants and parameter choices.
They also assisted with English translation, editing, LaTeX formatting, and bibliography preparation.
AI-generated bibliographic entries were checked by the first author.
The authors take full responsibility for the contents of this paper.





\clearpage

\clearpage
\appendix

\section*{Appendix}

\addtocontents{toc}{\protect\setcounter{tocdepth}{-1}}

\section{Proofs of Section~\ref{sec:analytic}}\label{app:analytic}

Throughout this appendix, $(M,g)$ and $p$ satisfy the hypotheses
of Proposition~\ref{prop:green}; no contractibility assumption is made.
\subsection{Proof of Proposition~\ref{prop:green}}
\begin{proof}
We first construct a smooth positive source $f$ satisfying
\[
 c\rho^{-3}\le f\le C\rho^{-3}.
\]
Since $M$ is second-countable and $d_g$ induces its manifold topology,
$(M,d_g)$ is separable. Moreover,
\[
 |\rho(x)^{-3}-\rho(y)^{-3}|
 \le 3|\rho(x)-\rho(y)|
 \le 3d_g(x,y),
\]
so $\rho^{-3}$ is Lipschitz.

Fix $0<\varepsilon<1$.
Applying \cite[Theorem 1]{Azagra} to $\rho^{-3}$
with the positive continuous error $\varepsilon\rho^{-3}$,
we obtain $f\in C^\infty(M)$ such that
\[
 |f-\rho^{-3}|\le\varepsilon\rho^{-3}.
\]
Consequently,
\[
 0<(1-\varepsilon)\rho^{-3}
 \le f\le(1+\varepsilon)\rho^{-3}
\]
by setting $c=1-\varepsilon$ and $C=1+\varepsilon$.

Write $\theta=\operatorname{AVR}(g)>0$.
For every $x\in M$ and $t>0$, Bishop--Gromov comparison
\cite[Lemma~7.1.4]{Petersen} gives
\[
 \theta\omega_4t^4
 \le\Vol_g(B_t^g(x))
 \le\omega_4t^4.
\]

Let $H_g$ be the heat kernel obtained by Dirichlet exhaustion,
and fix $0<\eta<1$.
By \cite[Corollary 3.1]{LY}, there is $A>0$ such that
\[
 H_g(x,y,t)
 \le
 \frac{A}{\Vol_g(B_{\sqrt t}^g(x))}
 \exp\left(-\frac{d_g(x,y)^2}{(4+\eta)t}\right).
\]
Consequently, for $x\ne y$,
\[
\begin{aligned}
 G_g(x,y)
 &:=
 \int_0^\infty H_g(x,y,t)\,dt\\
 &\le
 \frac{A}{\theta\omega_4}
 \int_0^\infty
 t^{-2}\exp\left(-\frac{d_g(x,y)^2}{(4+\eta)t}\right)\,dt\\
 &=\frac{A(4+\eta)}{\theta\omega_4d_g(x,y)^2}
 <\infty.
\end{aligned}
\]
The heat-kernel representation therefore gives a positive symmetric
Green kernel for $-\Delta_g$, agreeing with the Dirichlet exhaustion
construction \cite[Section~5]{LY}: on each bounded exhaustion domain,
the Dirichlet Green kernel is the time integral of its Dirichlet heat
kernel, and monotone convergence identifies the two exhaustion limits.

Now \cite[Theorem 5.2]{LY}, together with the volume bounds, yields
\[
\begin{aligned}
 G_g(x,y)
 &\ge
 C_{10}\int_{d_g(x,y)^2}^{\infty}
 \frac{dt}{\Vol_g(B_{\sqrt t}^g(x))}
 \ge\frac{C_{10}}{\omega_4d_g(x,y)^2},\\
 G_g(x,y)
 &\le
 C_{11}\int_{d_g(x,y)^2}^{\infty}
 \frac{dt}{\Vol_g(B_{\sqrt t}^g(x))}
 \le\frac{C_{11}}{\theta\omega_4d_g(x,y)^2}.
\end{aligned}
\]
Thus
\[
 a\,d_g(x,y)^{-2}\le G_g(x,y)\le b\,d_g(x,y)^{-2},
 \qquad x\ne y,
\]
where $a=C_{10}/\omega_4$ and
$b=C_{11}/(\theta\omega_4)$.

Define the source averages
\[
 k_f(x,t)=
 \frac{1}{\Vol_g(B_t^g(x))}
 \int_{B_t^g(x)}f(y)\,d\Vol_g(y).
\]
Polar-coordinate volume comparison
\cite[Lemma~7.1.2]{Petersen} yields, for every $s>0$,
\[
 \int_{B_s^g(p)}f(y)\,d\Vol_g(y)
 \le4C\omega_4\int_0^s\frac{r^3}{(1+r)^3}\,dr
 \le4C\omega_4s.
\]
If $t\le\rho(x)/2$, then $\rho(y)\ge\rho(x)/2$ on $B_t^g(x)$.
If $t>\rho(x)/2$, then $B_t^g(x)\subset B_{3t}^g(p)$.
Consequently,
\[
 k_f(x,t)\le
 \begin{cases}
 8C\rho(x)^{-3},&0<t\le\rho(x)/2,\\[1mm]
 12C\theta^{-1}t^{-3},&t>\rho(x)/2.
 \end{cases}
\]
Integrating these two ranges gives
\begin{equation}\label{eq:green-source-averages}
\begin{aligned}
 \int_0^\infty t\,k_f(x,t)\,dt
 &\le\frac{C(1+24/\theta)}{\rho(x)},\\
 \int_0^\infty k_f(x,t)\,dt
 &\le\frac{C(4+24/\theta)}{\rho(x)^2}.
\end{aligned}
\end{equation}

Set
\[
 P_f(x)=\int_MG_g(x,y)f(y)\,d\Vol_g(y).
\]
By Tonelli's theorem and
$d_g(x,y)^{-2}=2\int_{d_g(x,y)}^\infty t^{-3}\,dt$,
\[
\begin{aligned}
 P_f(x)
 &\le2b\int_0^\infty
 t^{-3}\int_{B_t^g(x)}f(y)\,d\Vol_g(y)\,dt\\
 &\le2b\omega_4\int_0^\infty t\,k_f(x,t)\,dt
 \le\frac{U}{\rho(x)}
\end{aligned}
\]
for some $U>0$. In particular, $P_f$ is finite everywhere.

For the lower bound, on $B_{\rho(x)}^g(x)$ we have
$d_g(x,y)<\rho(x)$ and $\rho(y)<2\rho(x)$. Hence
\[
 P_f(x)
 \ge\frac{ac}{8\rho(x)^5}\Vol_g(B_{\rho(x)}^g(x))
 \ge\frac{ac\theta\omega_4}{8\rho(x)}
 =\frac{L}{\rho(x)}.
\]

The Dirichlet exhaustion kernel $G_g$ is the minimal positive
Green kernel; in particular, $M$ is nonparabolic.
The volume and Green bounds give, with
$\sigma=\max\{1,b\omega_4,(a\theta\omega_4)^{-1}\}$,
\[
 \sigma^{-1}
 \frac{d_g(x,y)^2}{\Vol_g(B_{d_g(x,y)}^g(x))}
 \le G_g(x,y)\le
 \sigma
 \frac{d_g(x,y)^2}{\Vol_g(B_{d_g(x,y)}^g(x))},
 \qquad x\ne y.
\]
Together with
$\int_0^\infty k_f(p,t)\,dt<\infty$
from \eqref{eq:green-source-averages}, these verify the hypotheses
of \cite[Theorem 1.1]{NST}.

Let $w$ be the solution constructed there, satisfying
$\Delta_gw=f$ and $w(p)=0$.
Since
\[
 \int_M
 |G_g(p,y)-G_g(x,y)|f(y)\,d\Vol_g(y)
 \le P_f(p)+P_f(x)<\infty,
\]
the Green representation preceding \cite[equation (1.13)]{NST} gives
\[
\begin{aligned}
 w(x)
 &=w(x)-w(p)\\
 &=\int_M
 \bigl(G_g(p,y)-G_g(x,y)\bigr)f(y)\,d\Vol_g(y)\\
 &=P_f(p)-P_f(x).
\end{aligned}
\]

We make the initial $C^2$ regularity explicit.
Choose a smooth bounded exhaustion $\Omega_j$ containing $p$, and
let $G_j$ be its positive Dirichlet Green kernels. The normalized
classical Dirichlet potentials
\[
 w_j(x)=\int_{\Omega_j}
 \bigl(G_j(p,y)-G_j(x,y)\bigr)f(y)\,d\Vol_g(y)
 \qquad(x\in\Omega_j)
\]
satisfy $\Delta_gw_j=f$. Since $0\le G_j\le G_g$,
\[
 |w_j(x)|\le P_f(p)+P_f(x)\le2U.
\]
Monotone convergence applied to the two nonnegative potentials,
using $G_j\uparrow G_g$, gives $w_j(x)\to w(x)$ at every fixed $x$.
On nested relatively compact coordinate balls, the coefficients of
$\Delta_g$ are smooth and uniformly elliptic, and $f$ is smooth.
The interior Schauder estimate \cite[Theorem~6.2]{GT} therefore gives
locally uniform $C^{2,\alpha}$ bounds for any fixed $0<\alpha<1$;
a subsequence converges in $C^2$ on smaller balls to the same
pointwise limit $w$. Thus $w\in C^2(M)$.

Thus $u=P_f(p)-w=P_f$ satisfies $-\Delta_gu=f$.
Since $f$ and $g$ are smooth, interior elliptic regularity
\cite[Theorem~6.17]{GT} gives $u\in C^\infty(M)$.
The preceding estimates yield
\[
 \frac{L}{\rho}\le u\le\frac{U}{\rho},
\]
so $u$ is positive.

Finally, \cite[Theorem~1.3(i)]{NST} gives a constant
$K=K(4,\sigma)>0$ such that
\[
 |du|_g(x)=|dw|_g(x)
 \le K\int_0^\infty k_f(x,t)\,dt
 \le\frac{KC(4+24/\theta)}{\rho(x)^2},
\]
where the last inequality follows from
\eqref{eq:green-source-averages}.
Taking $D=KC(4+24/\theta)$ proves \eqref{eq:newton}.
\end{proof}

\subsection{Proof of Proposition~\ref{prop:conformal}}
\label{app:conformal}

\begin{proof}
Retain the potential and constants from Proposition~\ref{prop:green}.
Set
\[
 m=\frac{c}{U^3},\qquad B=\frac{D^2}{L^4},\qquad
 \varphi=\log u,\qquad \gbar=e^{2\varphi}g.
\]
Then
\[
 -\Delta_gu=f\ge\frac{c}{\rho^3}\ge mu^3.
\]

We use the conformal identities recalled in
Section~\ref{sec:preliminaries}.
For $\widetilde g=e^{2\eta}g$ and a smooth function $h$, they give
\begin{equation}\label{eq:conformal-hessian}
 \Hess_{\widetilde g}h
 =\Hess_gh-d\eta\otimes dh-dh\otimes d\eta
  +\langle d\eta,dh\rangle_g\,g.
\end{equation}
In dimension four,
\begin{equation}\label{eq:conformal-ricci}
 \Ric_{\widetilde g}
 =\Ric_g-2\Hess_g\eta+2d\eta\otimes d\eta
  -(\Delta_g\eta+2|d\eta|_g^2)g,
\end{equation}
and
\begin{equation}\label{eq:conformal-scalar}
 \Scal_{\widetilde g}
 =e^{-2\eta}
  \bigl(\Scal_g-6\Delta_g\eta-6|d\eta|_g^2\bigr).
\end{equation}

For $\gbar=u^2g$ we take $\eta=\varphi=\log u$.
Differentiation gives
\begin{equation}\label{eq:log-potential-derivatives}
\begin{aligned}
 d\varphi&=u^{-1}du,\\
 \Hess_g\varphi&=u^{-1}\Hess_gu-u^{-2}du\otimes du,\\
 \Delta_g\varphi&=u^{-1}\Delta_gu-u^{-2}|du|_g^2
                 =-\frac fu-u^{-2}|du|_g^2,\\
 |d\varphi|_g^2&=u^{-2}|du|_g^2.
\end{aligned}
\end{equation}
Substituting $\eta=\varphi$ in
\eqref{eq:conformal-ricci} and taking $h=\varphi$ in
\eqref{eq:conformal-hessian}, we obtain
\[
\begin{aligned}
 \Ric_{\gbar}
 &=\Ric_g-2\Hess_g\varphi
   +2d\varphi\otimes d\varphi
   -(\Delta_g\varphi+2|d\varphi|_g^2)g\\
 &=\Ric_g-\frac2u\Hess_gu+\frac4{u^2}du\otimes du
   -\left(\frac{\Delta_gu}{u}+\frac{|du|_g^2}{u^2}\right)g,\\
 \Hess_{\gbar}\varphi
 &=\Hess_g\varphi-2d\varphi\otimes d\varphi
   +|d\varphi|_g^2g\\
 &=u^{-1}\Hess_gu-3u^{-2}du\otimes du+u^{-2}|du|_g^2g.
\end{aligned}
\]
Adding these identities cancels the Hessian terms:
\begin{equation}\label{eq:weightedidentity}
\begin{aligned}
 \Ric_{\gbar}+2\Hess_{\gbar}\varphi
 +d\varphi\otimes d\varphi
 &=\Ric_g-d\varphi\otimes d\varphi-\Delta_g\varphi\,g\\
 &=\Ric_g+\frac fu g
 +u^{-2}\bigl(|du|_g^2g-du\otimes du\bigr)\\
 &\ge mu^2g=m\gbar.
\end{aligned}
\end{equation}
The last inequality uses $\Ric_g\ge0$ and
\[
 \bigl(|du|_g^2g-du\otimes du\bigr)(V,V)
 =|du|_g^2|V|_g^2-du(V)^2\ge0.
\]

For scalar curvature, \eqref{eq:conformal-scalar} and
\eqref{eq:log-potential-derivatives} give
\begin{equation}\label{eq:scalarbar}
\begin{aligned}
 \Scal_{\gbar}
 &=u^{-2}\bigl(\Scal_g-6\Delta_g\varphi-6|d\varphi|_g^2\bigr)\\
 &=u^{-2}\left(\Scal_g-6\frac{\Delta_gu}{u}
       +6\frac{|du|_g^2}{u^2}-6\frac{|du|_g^2}{u^2}\right)\\
 &=u^{-3}\bigl(-6\Delta_gu+\Scal_g\,u\bigr)\\
 &=u^{-3}(6f+\Scal_g\,u)\ge6m.
\end{aligned}
\end{equation}
Moreover,
\[
 |d\varphi|_{\gbar}^2
 =e^{-2\varphi}|d\varphi|_g^2
 =u^{-4}|du|_g^2
 \le\frac{\rho^4}{L^4}\frac{D^2}{\rho^4}
 =B.
\]

For completeness, let $\alpha:[0,1]\to M$ be a piecewise smooth
curve from $p$ to $x$. Since $\rho$ is $1$-Lipschitz for $d_g$,
$\rho\circ\alpha$ is absolutely continuous and
\[
 |(\rho\circ\alpha)'|\le|\alpha'|_g
 \quad\text{almost everywhere}.
\]
Thus
\[
 \operatorname{Length}_{\gbar}(\alpha)
 \ge L\int_0^1
 \frac{|(\rho\circ\alpha)'|}{\rho\circ\alpha}\,dt
 \ge L\log\rho(x),
\]
and hence
\begin{equation}\label{eq:bar-radial-distance}
 d_{\gbar}(p,x)\ge L\log\rho(x).
\end{equation}
Every closed $\gbar$-ball about $p$ is therefore a closed subset
of a compact $g$-ball. Here we use that $g$ and $\gbar$ induce
the same topology. Hopf--Rinow proves completeness of $(M,\gbar)$
\cite[Theorem 5.7.1]{Petersen}.

Finally,
\[
 \varphi\le\log U,\qquad
 -\varphi\ge\log\rho-\log U.
\]
Thus $-\varphi$ is bounded below and has compact sublevels.
\end{proof}

\subsection{Proof of Theorem~\ref{thm:joint}}
\label{app:morse}

Throughout this subsection, retain the data of
Proposition~\ref{prop:conformal} and put
\[
\begin{gathered}
 \psi=\frac{\varphi}{3},\qquad
 q=u^{2/3}g=e^{2\psi}g=e^{-4\varphi/3}\gbar,\\
 \Phi=\log U,\qquad
 \delta=\frac m8e^{-\sqrt B}>0.
\end{gathered}
\]
For this second conformal change, take $\eta=\psi=\varphi/3$ in
\eqref{eq:conformal-hessian} and \eqref{eq:conformal-scalar}.
Since $d\psi=d\varphi/3$ and $\Delta_g\psi=\Delta_g\varphi/3$,
\begin{equation}\label{eq:q-conformal-data}
\begin{aligned}
 \Hess_qh
 &=\Hess_gh-\frac13\left(
 d\varphi\otimes dh+dh\otimes d\varphi
 -\langle d\varphi,dh\rangle_g g\right),\\
 \Scal_q
 &=u^{-2/3}\left(\Scal_g-2\Delta_g\varphi
                        -\frac23|d\varphi|_g^2\right).
\end{aligned}
\end{equation}

\begin{lemma}[Proper ray Busemann functions]
\label{lem:proper-busemann}
Let $\gamma:[0,\infty)\to M$ be a $\gbar$-unit-speed ray
with $\gamma(0)=p$. Then
\[
 b(x)=\lim_{t\to\infty}
 \bigl(t-d_{\gbar}(x,\gamma(t))\bigr)
\]
is $1$-Lipschitz for $d_{\gbar}$, bounded below,
and has compact sublevels.
\end{lemma}

\begin{proof}
We first prove a uniform annular diameter bound.

The smooth 4-dimensional case of
\cite[Theorem 4.8]{CMS} gives, for all sufficiently large $r$,
a path-connected set $C_r$ such that
\[
 B_{4r}^g(p)\setminus B_{3r}^g(p)
 \subset C_r
 \subset B_{5r}^g(p)\setminus B_{2r}^g(p).
\]
Write $\theta=\operatorname{AVR}(g)>0$.
If $z_1,\ldots,z_N\in C_r$ are $r/4$-separated, then the
balls $B_{r/8}^g(z_i)$ are disjoint and contained in
$B_{41r/8}^g(p)$. Volume comparison gives
\[
 N\theta\omega_4(r/8)^4
 \le\Vol_g(B_{41r/8}^g(p))
 \le\omega_4(41r/8)^4.
\]
Thus $N\le41^4/\theta$, independently of $r$.

Choose a maximal such separated set. Its $r/4$-balls cover $C_r$,
and the intersection graph of this cover relative to $C_r$
is connected.
For $a,b\in C_r$ with $d_g(a,b)<r$, a minimizing $g$-segment
between them stays outside $B_r^g(p)$. Hence
\[
 d_{\gbar}(a,b)
 \le\frac{U}{1+r}\,d_g(a,b)<U.
\]
A simple chain in the finite covering graph uses at most $N-1$ edges.
Including the two endpoint-to-center segments therefore gives
\[
 \operatorname{diam}_{\gbar}(C_r)\le D_0,
 \qquad D_0=\bigl(\lceil41^4/\theta\rceil+2\bigr)U,
\]
independently of $r$.
For every sufficiently large $R$, take $r=R/3$. Since metric balls
are open,
\[
 \{x:d_g(p,x)=R\}
 \subset B_{4r}^g(p)\setminus B_{3r}^g(p)\subset C_r.
\]
Thus every such distance sphere has $\gbar$-diameter at most $D_0$.

The limit defining $b$ exists because its defining terms are
nondecreasing and bounded above by $d_{\gbar}(p,x)$.
The triangle inequality also gives
\[
 |b(x)-b(y)|\le d_{\gbar}(x,y).
\]

The ray $\gamma$ escapes every compact set, so
$d_g(p,\gamma(t))\to\infty$.
For $x$ outside a fixed compact set, choose $s\ge0$ such that
\[
 d_g(p,\gamma(s))=d_g(p,x).
\]
Then
\[
 d_{\gbar}(x,\gamma(s))\le D_0,
 \qquad
 d_{\gbar}(p,x)\le s+D_0,
\]
and consequently
\[
 b(x)\ge s-D_0
 \ge d_{\gbar}(p,x)-2D_0
 \ge L\log\rho(x)-2D_0.
\]
Together with continuity on the remaining compact set, this proves
boundedness below and compactness of all sublevels.
\end{proof}

\begin{lemma}[Weighted lower supports]
\label{lem:weighted-busemann}
Let $b$ be as in Lemma~\ref{lem:proper-busemann}.
For every compact set $K\subset M$, there is $H_K>0$ such that
each $x\in K$ admits a smooth local lower support $h$ for $b$:
\[
 h\le b\text{ near }x,\qquad h(x)=b(x).
\]
At the contact point,
\begin{equation}\label{eq:busemann-support-data}
\begin{gathered}
 |dh|_{\gbar}=1,\qquad
 \Hess_{\gbar}h(\nabla_{\gbar}h,\nabla_{\gbar}h)=0,\\
 \Delta_{\gbar}h
 -2\langle d\varphi,dh\rangle_{\gbar}\ge\delta,\\
 |\Hess_{\gbar}h|_{\gbar}\le H_K.
\end{gathered}
\end{equation}
\end{lemma}

\begin{proof}
The asymptotic-ray construction uses only completeness, as can be
seen directly from the proof of \cite[Proposition~7.3.8]{Petersen}.
Choose $t_j\to\infty$ and minimizing $\gbar$-segments
$\sigma_j:[0,a_j]\to M$ from $x$ to $\gamma(t_j)$, where
$a_j=d_{\gbar}(x,\gamma(t_j))\to\infty$.
A subsequence of the initial unit vectors converges, giving a
$\gbar$-unit-speed ray $\sigma$ from $x$.
For fixed $s\ge0$ and sufficiently large $j$,
\[
 t_j-d_{\gbar}(\sigma_j(s),\gamma(t_j))
 =t_j-a_j+s.
\]
Since the approximating Busemann functions are $1$-Lipschitz,
we can pass to the limit using $\sigma_j(s)\to\sigma(s)$ and obtain
\[
 b(\sigma(s))=b(x)+s.
\]
Our Busemann function is the negative of Petersen's convention;
no Ricci curvature assumption on $\gbar$ is used in this argument.
Set
\[
 \varphi_\sigma(s)=\varphi(\sigma(s)),\qquad
 T=\max\left\{
 2,\frac{32e^{\Phi-\varphi(x)+\sqrt B}}m
 \right\},
\]
and, for $0\le s\le T$, put
\[
 \ell(s)=1-\frac{s}{T},\qquad
 \zeta(s)=e^{(\varphi_\sigma(s)-\varphi_\sigma(0))/2}\ell(s).
\]
We have
\[
 \varphi_\sigma\le\Phi,\qquad |\varphi_\sigma'|\le\sqrt B,\qquad
 \zeta(0)=1,\qquad \zeta(T)=0.
\]

Let $\Pi_s:T_xM\to T_{\sigma(s)}M$ be parallel transport for
$\gbar$. For small $v\in T_xM$, define
\[
\begin{aligned}
 c_v(s)&=\exp_{\sigma(s)}^{\gbar}\bigl(\zeta(s)\Pi_sv\bigr),\\
 h(\exp_x^{\gbar}v)
 &=b(x)+T-\operatorname{Length}_{\gbar}(c_v).
\end{aligned}
\]
For sufficiently small $v$, the curves $c_v$ have nonzero velocity,
so their lengths depend smoothly on $v$.
Their initial and final points are $\exp_x^{\gbar}v$ and $\sigma(T)$.
Since $b(\sigma(T))=b(x)+T$ and $b$ is $1$-Lipschitz,
\[
 h(\exp_x^{\gbar}v)
 \le b(\sigma(T))
      -d_{\gbar}(\exp_x^{\gbar}v,\sigma(T))
 \le b(\exp_x^{\gbar}v).
\]
Equality holds at $v=0$.

Choose parallel $\gbar$-orthonormal fields
$E_1,E_2,E_3$ perpendicular to $\dot\sigma$ and put
\[
 I_i=
 \int_0^T
 \left(
  \zeta'^2
  -\zeta^2\,
   \gbar\bigl(R^{\gbar}(E_i,\dot\sigma)\dot\sigma,E_i\bigr)
 \right)\,ds.
\]
The first and second variation formulas
\cite[Lemma~5.4.2 and Theorem~6.1.4]{Petersen} give, at $x$,
\[
\begin{aligned}
 \nabla_{\gbar}h&=\dot\sigma(0),\\
 \Hess_{\gbar}h(\dot\sigma(0),\dot\sigma(0))&=0,\\
 \Delta_{\gbar}h&=-\sum_{i=1}^3I_i.
\end{aligned}
\]
Here the endpoint acceleration terms vanish: the initial endpoint
varies along a geodesic and the final endpoint is fixed.
For each perpendicular field $\zeta E_i$, the first variation of
speed vanishes pointwise, so the second variations of length and
of the energy $\frac12\int_0^T|\dot c|_{\gbar}^2\,ds$ agree.
For clarity, the radial Hessian identity follows directly from the
tangential variation. Extending $\sigma$ slightly backwards, we have
\[
 c_{\tau\dot\sigma(0)}(s)=\sigma(s+\tau\zeta(s)).
\]
For sufficiently small $|\tau|$, $1+\tau\zeta'(s)>0$, so
\[
 \operatorname{Length}_{\gbar}(c_{\tau\dot\sigma(0)})
 =\int_0^T(1+\tau\zeta'(s))\,ds=T-\tau.
\]
Hence $h(\exp_x^{\gbar}(\tau\dot\sigma(0)))=b(x)+\tau$.

Along $\sigma$, \eqref{eq:weightedidentity} gives
\[
 \Ric_{\gbar}(\dot\sigma,\dot\sigma)
 \ge m-2\varphi_\sigma''-(\varphi_\sigma')^2.
\]
Integration by parts therefore yields
\[
\begin{aligned}
 \sum_{i=1}^3I_i+2\varphi_\sigma'(0)
 &\le
 \int_0^T
 \left(
  3\zeta'^2-4\varphi_\sigma'\zeta\zeta'
  +(\varphi_\sigma')^2\zeta^2-m\zeta^2
 \right)\,ds\\
 &=
 \int_0^T e^{\varphi_\sigma-\varphi_\sigma(0)}
 \left(
  4\ell'^2
  -\left(\ell'+\frac{\varphi_\sigma'\ell}{2}\right)^2
  -m\ell^2
 \right)\,ds.
\end{aligned}
\]
Since $T\ge2$, on $[0,1]$ we have
\[
 \ell\ge\frac12,\qquad
 e^{\varphi_\sigma-\varphi_\sigma(0)}\ge e^{-\sqrt B}.
\]
Consequently,
\[
\begin{aligned}
 \sum_{i=1}^3I_i+2\varphi_\sigma'(0)
 &\le
 \frac{4e^{\Phi-\varphi_\sigma(0)}}T
 -\frac m4e^{-\sqrt B}\\
 &\le-\frac m8e^{-\sqrt B}
 =-\delta.
\end{aligned}
\]
Thus
\[
 \Delta_{\gbar}h
 -2\langle d\varphi,dh\rangle_{\gbar}
 =-\sum_{i=1}^3I_i-2\varphi_\sigma'(0)\ge\delta.
\]

For $x\in K$, the chosen lengths $T$ are uniformly bounded.
All corresponding coray segments lie in a compact
$\gbar$-neighborhood of $K$. On this neighborhood the curvature
is bounded, while $\zeta$ and $\zeta'$ are uniformly bounded.
More explicitly, let $T_K$ bound $T$, let $Z_0,Z_1$ bound
$|\zeta|,|\zeta'|$, and let $K_0$ bound the absolute sectional
curvature on this compact neighborhood. For $v\in T_xM$, write
$v^\perp=v-\gbar(v,\dot\sigma(0))\dot\sigma(0)$. The length second
variation depends only on this perpendicular part and gives
\[
\begin{aligned}
 \bigl|\Hess_{\gbar}h(v,v)\bigr|
 &\le\int_0^T(\zeta'^2+K_0\zeta^2)|v^\perp|_{\gbar}^2\,ds\\
 &\le T_K(Z_1^2+K_0Z_0^2)|v|_{\gbar}^2.
\end{aligned}
\]
In dimension four the tensor norm is at most twice the operator
norm, so $H_K=2T_K(Z_1^2+K_0Z_0^2)+1$ is an admissible bound in
\eqref{eq:busemann-support-data}.
\end{proof}

\begin{lemma}[Increasing reparametrization]
\label{lem:busemann}
There is a smooth strictly increasing function $\chi:\R\to\R$
such that $\chi\circ b$ is proper, bounded below, and belongs
to Shen's class $C(3)$ for $q$.
\end{lemma}

\begin{proof}
At each $x$, fix a calibrated coray and the lower support $h_x$
constructed in Lemma~\ref{lem:weighted-busemann}, and retain these
choices throughout.
Apply \eqref{eq:conformal-hessian} with base metric $\gbar$ and
conformal exponent $-2\varphi/3$, since $q=e^{-4\varphi/3}\gbar$.
This gives
\[
 \Hess_qh
 =
 \Hess_{\gbar}h
 +\frac23\left(
  d\varphi\otimes dh+dh\otimes d\varphi
  -\langle d\varphi,dh\rangle_{\gbar}\gbar
 \right).
\]
At a contact point, choose a $\gbar$-orthonormal frame
$e_0,e_1,e_2,e_3$ with $e_0=\nabla_{\gbar}h$, and write
\[
 A_{ij}=\Hess_qh(e_i,e_j).
\]
By \eqref{eq:busemann-support-data},
\begin{equation}\label{eq:busemann-tangential-trace}
 A_{11}+A_{22}+A_{33}
 =\Delta_{\gbar}h
  -2\langle d\varphi,dh\rangle_{\gbar}
 \ge\delta.
\end{equation}

Choose a lower bound $b_*$ for $b$ and set
\[
 K_n=\{x:b(x)\le b_*+n+2\},
 \qquad n=0,1,\ldots.
\]
These sets are compact. Lemma~\ref{lem:weighted-busemann}
and the gradient bound for $\varphi$ give constants $C_n\ge1$
such that
\[
 |A_{ij}|\le C_n\qquad(x\in K_n).
\]
Enlarge them so that $C_n$ is nondecreasing.

We use the following elementary matrix estimate.
If $A$ is symmetric, $|A_{ij}|\le C_n$, and
$A_{11}+A_{22}+A_{33}\ge\delta$, then
\begin{equation}\label{eq:busemann-matrix-correction}
 \lambda\ge4C_n+\frac{6C_n^2}{\delta}
 \quad\Longrightarrow\quad
 \tr_P^{\gbar}\bigl(\Hess_qh+\lambda\,dh\otimes dh\bigr)
 \ge\frac{\delta}{2}
\end{equation}
for every three-plane $P\subset T_xM$.
Indeed, put $S=A_{11}+A_{22}+A_{33}$ and let
$\nu=(\nu_0,\nu')$ be a $\gbar$-unit normal to $P$.
Cauchy--Schwarz and Young's inequality give
\[
\begin{aligned}
 \sum_{i,j=1}^3\nu_iA_{ij}\nu_j
 &\le3C_n|\nu'|^2,\\
 2\left|\nu_0\sum_{i=1}^3A_{0i}\nu_i\right|
 &\le\frac{\delta}{2}\nu_0^2
     +\frac{6C_n^2}{\delta}|\nu'|^2.
\end{aligned}
\]
Since $A_{00}\ge-C_n$,
\[
\begin{aligned}
 &\tr_P^{\gbar}\bigl(\Hess_qh+\lambda\,dh\otimes dh\bigr)\\
 &\qquad=
 S\nu_0^2+(S+A_{00}+\lambda)|\nu'|^2
 -\sum_{i,j=1}^3\nu_iA_{ij}\nu_j
 -2\nu_0\sum_{i=1}^3A_{0i}\nu_i\\
 &\qquad\ge
 \frac{\delta}{2}\nu_0^2
 +\left(\delta+\lambda-4C_n
              -\frac{6C_n^2}{\delta}\right)|\nu'|^2\\
 &\qquad\ge\frac{\delta}{2}.
\end{aligned}
\]

Choose a smooth positive function $\Lambda:\R\to(0,\infty)$
such that, for every $n=0,1,\ldots$,
\[
 \Lambda(t)\ge4C_n+\frac{6C_n^2}{\delta}
 \qquad\text{for }b_*+n\le t\le b_*+n+1.
\]
Such a function exists by taking a smooth majorant of these
locally bounded stripwise constants.
Define
\[
\begin{aligned}
 \chi'(t)&=\exp\left(\int_{b_*}^t\Lambda(s)\,ds\right),\\
 \chi(t)&=b_*+\int_{b_*}^t\chi'(s)\,ds.
\end{aligned}
\]
Then $\chi'>0$ and $\chi''/\chi'=\Lambda$ on $\R$.
For every $t\ge b_*$, we also have
\[
 \chi'(t)\ge1,\qquad\chi(t)\ge t.
\]
In particular, $\chi\circ b$ is bounded below and has compact
sublevels.

Since $\chi$ is increasing, $\chi\circ h$ is a smooth local
lower support for $\chi\circ b$.
At its contact point,
\[
 \Hess_q(\chi\circ h)
 =
 \chi'(b(x))
 \left(
  \Hess_qh+
  \frac{\chi''(b(x))}{\chi'(b(x))}\,dh\otimes dh
 \right).
\]
The $q$-orthonormal frame corresponding to $e_i$ is
$e^{2\varphi(x)/3}e_i$. Thus
\eqref{eq:busemann-matrix-correction} gives
\[
 \tr_P^q\Hess_q(\chi\circ h)
 \ge
 \frac{\delta}{2}
 e^{4\varphi(x)/3}\chi'(b(x))>0
 \qquad(\dim P=3).
\]

Fix $x_0\in M$ and choose an open neighborhood $W$ with
$\overline W\subset\intt K_N$ for some $N$.
At each $x\in W$, choose the lower support $h_x$ above and put
$S_x=\Hess_q(\chi\circ h_x)_x$.
The three-plane estimate and the coefficient bound give
\[
 \tr_P^q S_x\ge\eta_0,\qquad
 |S_x(v,v)|\le A_0|v|_q^2,
\]
for every $x\in W$, every three-plane $P\subset T_xM$, and every
$v\in T_xM$. The constants are uniform on $W$; one may take
\[
\begin{aligned}
 \eta_0&=\min_{K_N}
 \frac{\delta}{2}e^{4\varphi/3}\chi'(b)>0,\\
 A_0&=1+\max_{K_N}
 e^{4\varphi/3}\chi'(b)\bigl(4C_N+\Lambda(b)\bigr).
\end{aligned}
\]
Indeed, $|A_{ij}|\le C_N$ implies
$|\Hess_qh_x(v,v)|\le4C_Ne^{4\varphi(x)/3}|v|_q^2$,
and $|dh_x|_q^2=e^{4\varphi(x)/3}$.

We spell out the perturbation argument of
\cite[Section 3.1, Lemma 6]{Shen}.
Choose
\[
 0<\varepsilon\le
 \min\left\{\frac16,\frac{\eta_0}{18A_0}\right\}.
\]
For any triple $v_1,v_2,v_3\in T_xM$ with
$|q(v_i,v_j)-\delta_{ij}|<\varepsilon$, its Gram matrix $G$
satisfies $\|G-I\|_{\mathrm{op}}\le3\varepsilon\le1/2$.
The vectors
\[
 w_i=\sum_{j=1}^3v_j(G^{-1/2})_{ji}
\]
are $q$-orthonormal. Setting $T=(S_x(w_i,w_j))$, we obtain
\[
\begin{aligned}
 \sum_{i=1}^3S_x(v_i,v_i)
 &=\tr(TG)\\
 &\ge\tr T-3\|T\|_{\mathrm{op}}\|G-I\|_{\mathrm{op}}\\
 &\ge\eta_0-9A_0\varepsilon
 \ge\frac{\eta_0}{2}.
\end{aligned}
\]

Taylor expansion of the smooth lower support along $q$-geodesics gives
\[
\begin{aligned}
 C_{\chi\circ b}^q(x;v)
 &:={\liminf}_{t\to0}
 \frac{
  (\chi\circ b)(\exp_x^q(tv))
  +(\chi\circ b)(\exp_x^q(-tv))
  -2(\chi\circ b)(x)
 }{t^2}\\
 &\ge S_x(v,v).
\end{aligned}
\]
Consequently, for every $x\in W$ and every such triple,
\[
 \sum_{i=1}^3C_{\chi\circ b}^q(x;v_i)
 \ge\frac{\eta_0}{2}>\frac{\eta_0}{4}.
\]
The function $b$ is $1$-Lipschitz for $\gbar$.
Local comparison of $q$ and $\gbar$ and smoothness of $\chi$
therefore make $\chi\circ b$ locally Lipschitz for $q$.
These are the locally uniform conditions defining $C(3)$ in
\cite[Section 3.1]{Shen}. Hence $\chi\circ b\in C(3)$ on $M$.
\end{proof}

\begin{proof}[Proof of Theorem~\ref{thm:joint}]
Use the same potential $u$ as in Proposition~\ref{prop:conformal}.
Since $u\le U$,
\[
 q=u^{-4/3}\gbar\ge U^{-4/3}\gbar.
\]
Consequently,
\[
 d_q(p,x)\ge U^{-2/3}d_{\gbar}(p,x).
\]
A closed $q$-ball is contained in a compact closed $\gbar$-ball
and is closed in their common manifold topology. It is therefore
compact, and Hopf--Rinow gives completeness of $q$.

Since $(M,\gbar)$ is complete and noncompact, choose a
$\gbar$-unit-speed ray $\gamma$ issuing from $p$, and let $b$
be its Busemann function.

By Lemma~\ref{lem:proper-busemann} and Lemma~\ref{lem:busemann},
the function $\chi\circ b$ is proper, bounded below, and belongs
to $C(3)$ for $q$.
Apply \cite[Section 3.1, Lemma 5]{Shen} with approximation
error equal to $1$. We obtain a smooth proper Morse function
$F$ such that
\[
 |F-\chi\circ b|<1,\qquad F\in C(3).
\]
The error bound makes $F$ bounded below, so $F$ is an exhaustion.
For a smooth function, Shen's $C(3)$ condition is equivalent to
\[
 \lambda_{(3)}^q(\Hess_qF)>0.
\]
At a critical point, three negative Hessian eigenvalues would
contradict this inequality. Hence
\[
 \ind_xF\le2\qquad(x\in\Crit F).
\]
\end{proof}

\section{Proofs of Section~\ref{sec:finite-exteriors}}
\label{app:finite-exteriors}

\subsection{Proof of Theorem~\ref{thm:cores}}

\begin{lemma}[Intersection of fundamental-group images]
\label{lem:pi1-image-intersection}
Let $X=U\cup V$, where $U$ and $V$ are open subsets of $X$
such that $U$, $V$, and $U\cap V$ are path connected.
For a common basepoint $z\in U\cap V$,
\[
\begin{aligned}
 &\operatorname{im}\bigl(\pi_1(U,z)\to\pi_1(X,z)\bigr)
 \cap
 \operatorname{im}\bigl(\pi_1(V,z)\to\pi_1(X,z)\bigr)\\
 &\qquad=
 \operatorname{im}\bigl(\pi_1(U\cap V,z)\to\pi_1(X,z)\bigr).
\end{aligned}
\]
\end{lemma}

\begin{proof}
Van Kampen identifies $\pi_1(X,z)$ with the group pushout
\cite[Theorem~1.20]{Hatcher}.
Write
\[
 A=\pi_1(U,z),\qquad
 B=\pi_1(V,z),\qquad
 C=\pi_1(U\cap V,z),
\]
and denote the pushout by $P$.
Let $\overline A,\overline B,\overline C$ be their actual images
in $P$, and put
\[
 P'=\overline A*_{\overline C}\overline B.
\]
The inclusions of these images give a homomorphism $P'\to P$.
Conversely, the maps
$A\to\overline A\to P'$ and $B\to\overline B\to P'$
agree on $C$, so the original universal property gives $P\to P'$.
The composites are the identity on their generating vertex groups,
hence are inverse isomorphisms.
The edge maps defining $P'$ are injective.
The desired identity therefore follows from the normal-form theorem
for amalgamated free products
\cite[Chapter~IV, Theorem~2.6, p.~187]{LyndonSchupp}:
if $a\in\overline A\setminus\overline C$ and
$b\in\overline B\setminus\overline C$ represented the same element,
then $ab^{-1}$ would be a reduced word of length two representing
the identity, contrary to that theorem.
\end{proof}

\begin{proof}[Proof of Theorem~\ref{thm:cores}]
Let $H$ be a smooth strictly three-convex Morse exhaustion for $q$.
Three negative critical Hessian eigenvalues would contradict
$\lambda_{(3)}^q(\Hess_qH)>0$, so every critical index is at most two.
Choose strictly increasing regular values $t_i\to\infty$
with $H(p)<t_1$, and put
\[
 K_i=\{H\le t_i\},\qquad Z_i=M\setminus K_i.
\]
Let $D_i$ be the component of $K_i$ containing $p$.
Each $K_i$ is a compact regular domain, possibly disconnected,
and each $D_i$ is a connected compact regular domain.

Since $t_i<t_{i+1}$,
\[
 D_i\subset\intt D_{i+1}.
\]
For any $x\in M$, choose a path from $p$ to $x$.
Its image is compact, so for sufficiently large $i$ it lies
in $\{H<t_i\}$. Hence
\[
 \bigcup_{i\ge1}\intt D_i=M.
\]

The Morse handle decomposition recalled in
Section~\ref{sec:preliminaries} builds each $K_i$ from
finitely many handles of indices $0,1,2$.
Each connected component has connected boundary:
$0$-handles start with $S^3$, $1$-handles merge or preserve
connected boundaries, and $2$-handle surgery along a circle
does not disconnect a connected three-manifold.
Thus $Y_{D_i}$ is connected.

Removing finitely many disjoint compact domains with connected
boundaries leaves a path-connected complement.
Indeed, choose a path between two exterior points transverse
to their boundaries.
It has finitely many boundary crossings; replace each excursion
through a removed domain by a path in its connected boundary,
pushed slightly into the exterior collar.
This proves that $Z_i$ is path connected.
Filling the components of $K_i$ other than $D_i$ shows that
$O_{D_i}$ is path connected as well.

Choose a minimizing $g$-ray $\eta$ from $p$ and strictly
increasing times $\tau_i$ such that
\[
 K_i\subset B_{\tau_i}^g(p).
\]
Set $x_i=\eta(\tau_i)$. Then
\[
 \eta([\tau_i,\infty))\subset Z_i\subset O_{D_i}.
\]
We use these points and the corresponding ray segments
for all basepoint changes below.

Reading the Morse cobordism between $t_i$ and $t_j$ in reverse
replaces an index-$\lambda$ handle by an index-$(4-\lambda)$ handle.
All these dual indices are at least two, so attaching them creates
no fundamental-group generators.
Using boundary collars and the chosen basepoint paths, we obtain
\[
 \pi_1(Z_j,x_j)\twoheadrightarrow\pi_1(Z_i,x_i)
 \qquad(j\ge i).
\]
This is the reverse-handle argument in
\cite[proof of Proposition~4.3]{CG}.

Similarly, turning each component of $K_i$ other than $D_i$
upside down gives a boundary-relative handle decomposition
with indices $2,3,4$.
Filling those components therefore gives
\[
 \pi_1(Z_i,x_i)\twoheadrightarrow\pi_1(O_{D_i},x_i).
\]
Since $Z_j\subset O_{D_j}$, the composite surjection factors as
\[
 \pi_1(Z_j,x_j)\longrightarrow
 \pi_1(O_{D_j},x_j)\longrightarrow
 \pi_1(O_{D_i},x_i).
\]
Hence
\begin{equation}\label{eq:cores-exterior-surjection}
 \pi_1(O_{D_j},x_j)\twoheadrightarrow\pi_1(O_{D_i},x_i)
 \qquad(j\ge i).
\end{equation}

In the volume-growth normalization of \cite{CMS}, the parameter is
$v=\theta\omega_4$, since
$\lim_{R\to\infty}\Vol_g(B_R^g(p))/R^4=\theta\omega_4$.
The smooth 4-dimensional case of \cite[Theorem~4.8]{CMS}
bounds the image of the fundamental group of the annular component
$C_r\subset B_{5r}^g(p)\setminus B_{2r}^g(p)$ containing
$B_{4r}^g(p)\setminus B_{3r}^g(p)$.
As in the beginning of the proof of
\cite[Proposition~6.6]{CMS}, choose an open path-connected
$A_r\subset C_r$ containing that smaller annulus.
Thus, for every sufficiently large $r$, we have
\[
 B_{4r}^g(p)\setminus B_{3r}^g(p)
 \subset A_r
 \subset B_{5r}^g(p)\setminus B_{2r}^g(p)
\]
and, by factoring the inclusion through $C_r$,
\begin{equation}\label{eq:cores-annular-image}
 \left|
 \operatorname{im}\left[
  \pi_1(A_r,z)\longrightarrow
  \pi_1(B_{6r}^g(p)\setminus B_r^g(p),z)
 \right]\right|\le N,
 \qquad(z\in A_r),
\end{equation}
where $N=\lceil C(4,v)\rceil$ is independent of $r$ and of $H$.
Fix $R_0>0$ so that this construction applies for every $r>R_0$.

Fix $i$ and an arbitrary finite subset
$S\subset\pi_1(O_{D_i},x_i)$.
Choose based loops representing its elements, and take $r>R_0$
large enough that
\[
 D_i\subset B_r^g(p),\qquad \tau_i<3r,
\]
and all these loops lie in $B_{3r}^g(p)$.
Set $z=\eta(7r/2)$ and define
\[
\begin{aligned}
 \mathcal U_r&=(A_r\cup B_{4r}^g(p))\setminus D_i,\\
 \mathcal V_r&=A_r\cup
              (M\setminus\overline{B_{3r}^g(p)}).
\end{aligned}
\]
These are open subsets of $O_{D_i}$ satisfying
\[
 \mathcal U_r\cup\mathcal V_r=O_{D_i},
 \qquad
 \mathcal U_r\cap\mathcal V_r=A_r.
\]

Both sets are path connected.
Indeed, every point of $\mathcal U_r\setminus A_r$ has
$g$-distance less than $3r$ from $p$, whereas every point of
$\mathcal V_r\setminus A_r$ has distance at least $4r$.
Join such a point to $z$ by a path in $O_{D_i}$ and stop at
its first intersection with the sphere of radius $7r/2$.
In the first case the initial segment lies in
$B_{4r}^g(p)\setminus D_i$; in the second case it lies in
$M\setminus\overline{B_{3r}^g(p)}$.
The stopping point belongs to $A_r$ and can be joined to $z$
inside $A_r$.

Transport $S$ to $\pi_1(O_{D_i},z)$ along the ray segment.
Its chosen representatives, including their basepoint paths,
lie in $\mathcal U_r$.

Choose $j\ge i$ such that
\[
 \overline{B_{4r}^g(p)}\subset\intt D_j.
\]
Since $D_j\subset K_j\subset B_{\tau_j}^g(p)$,
we have $\tau_j>4r$.
By \eqref{eq:cores-exterior-surjection}, the same transported
classes have representatives in $O_{D_j}$, joined to $z$
by the ray segment from $z$ to $x_j$.
Since $O_{D_j}$ lies outside $\overline{B_{4r}^g(p)}$
and this ray segment has radial coordinate at least $7r/2$,
the resulting representatives lie in $\mathcal V_r$.

Lemma~\ref{lem:pi1-image-intersection} therefore places every
transported element of $S$ in
\[
 \operatorname{im}\bigl(
 \pi_1(A_r,z)\to\pi_1(O_{D_i},z)
 \bigr).
\]
The latter map factors through
\[
 A_r\subset B_{6r}^g(p)\setminus B_r^g(p)\subset O_{D_i}.
\]
Thus \eqref{eq:cores-annular-image} gives $|S|\le N$.
Since $S$ was an arbitrary finite subset,
\[
 |\pi_1(O_{D_i},x_i)|\le N.
\]
The same $N$ works for every $i$, so $N_H=N$ can be chosen
to depend only on $(M,g)$, independently of $H$.

Finally, $Y_{D_i}\subset H^{-1}(t_i)$ is a regular level
with outward normal $\nabla_qH/|dH|_q$.
Hence
\[
 H_q(Y_{D_i})
 =
 \frac{\tr_{(\nabla_qH)^\perp}^q\Hess_qH}{|dH|_q}
 \ge
 \frac{\lambda_{(3)}^q(\Hess_qH)}{|dH|_q}
 >0.
\]
\end{proof}

\section{Proofs of Section~\ref{sec:cuts}}\label{app:union}

Throughout this section, retain the conformal and Morse data
and the compact regular domains $C_i$ constructed above.
In particular, $q$ and $\gbar$ are complete, and $F$ is
a smooth strictly three-convex Morse exhaustion for $q$.

\subsection{Proof of Proposition~\ref{prop:forcing}}

\begin{proof}
Choose $\sigma\in C^\infty(\R)$ such that
\[
 0\le\sigma\le1,\qquad
 \sigma(t)=0\quad(t\le1/3),\qquad
 \sigma(t)=1\quad(t\ge2/3),
\]
and put $C_\sigma=\sup_\R|\sigma'|$.
Choose
\[
 R>\max\left\{6,\frac{2\pi C_\sigma}{\sqrt m}\right\}.
\]
This choice is independent of the plateau values below.

Set $A=C_i$ and
\[
 r_A(x)=\inf_{y\in A}d_{\gbar}(x,y).
\]
Completeness of $\gbar$ and compactness of $A$ imply that
$\{r_A\le R+2\}$ is compact.
Choose $j>i$ such that
\[
 \{r_A\le R+2\}\subset\intt C_j,
\]
and put $B=C_j$.
Thus $r_A>R+2$ on $Y_B$.

The distance function $r_A$ is $1$-Lipschitz.
The Riemannian metric $\gbar$ induces the second-countable manifold
topology, so $(M,d_{\gbar})$ is separable.
Apply \cite[Theorem~1]{Azagra} with constant value error $1$
and Lipschitz-loss allowance $1$. There is $r\in C^\infty(M)$ with
\[
 |r-r_A|\le1,\qquad |dr|_{\gbar}\le2.
\]
Consequently,
\[
 r\le1<R/3\quad\text{on }A,\qquad
 r>R+1>2R/3\quad\text{on }Y_B.
\]

Choose constants $a_0,a_1$ such that
\[
\begin{gathered}
 a_0>\max_{Y_A}u^{-2/3}H_q(Y_A),\\
 0<a_1<
 \min\left\{a_0,\min_{Y_B}u^{-2/3}H_q(Y_B)\right\}.
\end{gathered}
\]
The second minimum is positive by Theorem~\ref{thm:cores}.
Put
\[
\begin{aligned}
 \beta_\ell&=\arctan\frac{a_\ell}{\sqrt m}\quad(\ell=0,1),\\
 \beta(t)&=(1-\sigma(t/R))\beta_0+\sigma(t/R)\beta_1,
\end{aligned}
\]
and define
\[
 a(t)=\sqrt m\,\tan\beta(t),\qquad h=a\circ r.
\]
Since
\[
 0<\beta_1\le\beta(t)\le\beta_0<\pi/2,
\]
both functions are smooth and positive.
Moreover,
\[
 a=a_0\quad(t\le R/3),\qquad
 a=a_1\quad(t\ge2R/3).
\]
Hence $h=a_0$ on $A$ and $h=a_1$ on $Y_B$, giving
\[
 H_q(Y_A)<u^{2/3}h,\qquad
 H_q(Y_B)>u^{2/3}h.
\]

Since $|\beta_0-\beta_1|<\pi/2$,
\[
 |\beta'|
 \le\frac{\pi C_\sigma}{2R}
 <\frac{\sqrt m}{4}.
\]
Thus
\[
 a'=\frac{a^2+m}{\sqrt m}\,\beta',
 \qquad |a'|\le\frac{a^2+m}{4}.
\]
Using $|dr|_{\gbar}\le2$, we obtain
\[
 2|dh|_{\gbar}
 =2|a'(r)|\,|dr|_{\gbar}
 \le h^2+m.
\]
This proves \eqref{eq:forcing}.

To prove \eqref{eq:DEC}, put
\[
 z=u^{-4/3}|du|_g,\qquad k=u^{2/3}h.
\]
The conformal scalar-curvature formula from
Section~\ref{sec:preliminaries}, together with $-\Delta_gu=f$,
gives
\[
\begin{aligned}
 \Scal_q
 &=u^{-2/3}\Scal_g+2fu^{-5/3}
   +\frac43u^{-8/3}|du|_g^2\\
 &\ge2mu^{4/3}+\frac43z^2.
\end{aligned}
\]
Also,
\[
\begin{aligned}
 |dk|_q
 &\le\frac23hu^{-2/3}|du|_g+u^{1/3}|dh|_g\\
 &=\frac23kz+u^{4/3}|dh|_{\gbar}.
\end{aligned}
\]
Therefore,
\[
\begin{aligned}
 \Scal_q+\frac43k^2-2|dk|_q
 &\ge
 \frac43\left(z-\frac k2\right)^2
 +u^{4/3}\bigl(2m+h^2-2|dh|_{\gbar}\bigr)\\
 &\ge mu^{4/3}>0.
\end{aligned}
\]
\end{proof}

\subsection{Proof of Theorem~\ref{input:eichmair}}

\begin{proof}
Fix the given orientation of $M$.
In our applications, an orientation exists because $M$ is contractible.
Put $\mathsf p=(\psi/3)q$.
For every hypersurface $\Sigma\subset B$,
\[
 \tr_{T\Sigma}^q\mathsf p=\psi>0,
 \qquad
 |\tr_{T\Sigma}^q\mathsf p|=\psi.
\]
Apply \cite[Theorem~5.1]{EichmairHorizons} to the complete
initial data set $(M,q,\mathsf p)$, with $n=4$ and
\[
 \Omega=\intt B\setminus A,\qquad
 \partial_1\Omega=Y_B,\qquad
 \partial_2\Omega=Y_A.
\]
The normal pointing out of $\Omega$ along $Y_B$ is the outward
normal of $B$, and the normal pointing into $\Omega$ along $Y_A$
is the outward normal of $A$.
The hypotheses therefore read
\[
 H_q(Y_B)>\psi,\qquad H_q(Y_A)\le\psi.
\]

The theorem gives an open set $\Omega'\subset\Omega$ with
\[
 \partial\Omega'=Y_B\sqcup\Sigma',
\]
where $\Sigma'$ is an embedded $C^{2,\alpha}$ hypersurface
satisfying $H_q=\psi$ with respect to the normal pointing into
$\Omega'$.
Since $\psi>0$, the absolute value in the equation is simply
$|\tr_{T\Sigma'}^q\mathsf p|=\psi$, with no loss of regularity.
In local graph coordinates, the $C^{2,\alpha}$ graph function
satisfies a quasilinear uniformly elliptic equation on smaller
neighborhoods. Evaluating its coefficients along this graph gives
$C^{1,\alpha}$ coefficients and right-hand side.
The linear interior regularity theorem
\cite[Theorem~6.17]{GT} first gives $C^{3,\alpha}$ regularity;
iteration gives smoothness because $q$ and $\psi$ are smooth.

Put
\[
 E=\intt B\setminus\Omega'.
\]
Near $Y_B$, the smooth domain $\Omega'$ occupies the inward side
of $Y_B$.
Compactness gives an inward collar of $Y_B$ contained in $\Omega'$.
Thus $E$ stays away from $Y_B$ and is compact in $M$.
Moreover,
\[
 A\subseteq E\subset\intt B,\qquad Y_E=\Sigma',
\]
and the outward normal of $E$ points into $\Omega'$.
Hence $E$ is a compact regular domain with $H_q(Y_E)=\psi$.
Its boundary is nonempty because
$\varnothing\ne A\subseteq E\ne M$ and $M$ is connected.

Let $D$ be a compact regular domain with
\[
 A\subseteq D\subset\intt B,\qquad H_q(Y_D)\le\psi.
\]
Then
\[
 \Omega''=\intt B\setminus D\subset\Omega
\]
has boundary $Y_B\sqcup Y_D$, and the normal pointing into
$\Omega''$ is the outward normal of $D$.
The outermost property in \cite[Theorem~5.1]{EichmairHorizons}
gives $\Omega'\subset\Omega''$, hence $D\subseteq E$.

Finally, compactness of $Y_A$ and the strict inequality
$H_q(Y_A)<\psi$ give a sufficiently small $t>0$ such that
\[
 A_t=\{x:d_q(x,A)\le t\}
\]
is a compact regular domain in $\intt B$ with
$H_q(Y_{A_t})<\psi$.
The enclosure property gives $A_t\subseteq E$, and hence
$A\subset\intt E$.
\end{proof}

\subsection{Proof of Proposition~\ref{prop:outerpsc}}

\begin{proof}
The function $k=u^{2/3}h$ is positive, and \eqref{eq:forcing} gives
\[
 H_q(Y_A)<k,\qquad H_q(Y_B)>k.
\]
Theorem~\ref{input:eichmair}, applied with $\psi=k$, gives
a compact regular domain $T$ with
\[
 A\subset\intt T,\qquad T\subset\intt B,\qquad H_q(Y_T)=k,
\]
containing every compact regular domain $L$ with
\[
 A\subseteq L\subset\intt B,\qquad H_q(Y_L)\le k.
\]
If moreover $T\subseteq L$, then $L\subseteq T$, so $L=T$.
Thus $T$ is weakly outermost in the sense of
Definition~\ref{def:outermost}.

We verify PSC on each connected component $\Sigma$ of $Y_T$.
Write $q_\Sigma=q|_{T\Sigma}$ and let $\nu$ be the outward
$q$-unit normal of $T$.
Under a normal variation with speed $\eta$,
\[
\begin{aligned}
 H_q'&=-\Delta_{q_\Sigma}\eta
       -\bigl(\Ric_q(\nu,\nu)+|\mathrm{II}_q|_{q_\Sigma}^2\bigr)\eta,\\
 k'&=\eta\,\nu k.
\end{aligned}
\]
Consequently, the linearization of $H_q-k$ is
\[
 \mathcal L_\Sigma\eta
 =-\Delta_{q_\Sigma}\eta
  -\bigl(\Ric_q(\nu,\nu)+|\mathrm{II}_q|_{q_\Sigma}^2
         +\nu k\bigr)\eta.
\]

Its first eigenvalue is nonnegative.
Indeed, a negative first eigenvalue and its positive eigenfunction
would give a small outward deformation of $\Sigma$ with
$H_q-k<0$.
Leaving the other boundary components fixed would produce a larger
admissible compact regular domain, contradicting weak outermostness.
This is the principal-eigenfunction variation argument in
\cite[proof of Theorem~2.1]{GallowaySchoen}, applied here directly
to $H_q-k$ and the self-adjoint operator $\mathcal L_\Sigma$.
The PSC conclusion below follows from our displayed stability
inequality and \eqref{eq:DEC}.

It follows that, for every $\eta\in C^\infty(\Sigma)$,
\[
 \int_\Sigma
 \left(
 |d\eta|_{q_\Sigma}^2
 -\bigl(\Ric_q(\nu,\nu)+|\mathrm{II}_q|_{q_\Sigma}^2
        +\nu k\bigr)\eta^2
 \right)d\Vol_{q_\Sigma}\ge0.
\]
The Gauss equation and $H_q(\Sigma)=k$ give
\[
 2\Ric_q(\nu,\nu)
 =\Scal_q-\Scal_{q_\Sigma}
   +k^2-|\mathrm{II}_q|_{q_\Sigma}^2.
\]
Using
\[
 |\mathrm{II}_q|_{q_\Sigma}^2\ge k^2/3,
 \qquad \nu k\ge-|dk|_q,
\]
we obtain
\[
\begin{aligned}
 &\int_\Sigma
 \left(2|d\eta|_{q_\Sigma}^2+\Scal_{q_\Sigma}\eta^2\right)
 d\Vol_{q_\Sigma}\\
 &\qquad\ge
 \int_\Sigma
 \left(\Scal_q+\frac43k^2-2|dk|_q\right)\eta^2
 d\Vol_{q_\Sigma}.
\end{aligned}
\]
By compactness and \eqref{eq:DEC},
\[
 \delta_\Sigma
 =\min_\Sigma\left(\Scal_q+\frac43k^2-2|dk|_q\right)>0.
\]
Adding the nonnegative term
$6\int_\Sigma|d\eta|_{q_\Sigma}^2\,d\Vol_{q_\Sigma}$ to the
left-hand side gives
\[
 \int_\Sigma
 \left(8|d\eta|_{q_\Sigma}^2+\Scal_{q_\Sigma}\eta^2\right)
 d\Vol_{q_\Sigma}
 \ge\delta_\Sigma\int_\Sigma\eta^2\,d\Vol_{q_\Sigma}.
\]
Thus the conformal Laplacian
\[
 -8\Delta_{q_\Sigma}+\Scal_{q_\Sigma}
\]
has first eigenvalue $\lambda_\Sigma\ge\delta_\Sigma>0$.
Let $\zeta>0$ be a corresponding smooth eigenfunction.
The three-dimensional conformal scalar-curvature formula yields
\[
 \Scal_{\zeta^4q_\Sigma}
 =\zeta^{-5}
   \bigl(-8\Delta_{q_\Sigma}\zeta+\Scal_{q_\Sigma}\zeta\bigr)
 =\lambda_\Sigma\zeta^{-4}>0.
\]
Therefore every component of $Y_T$ admits PSC.
\end{proof}

\subsection{Proof of Lemma~\ref{lem:psc-collar}}

\begin{lemma}[Inward collar data]\label{lem:collarcut}
Let $K$ be a connected compact regular domain with
$H_q(Y_K)>0$, and let $A\subset\intt K$ be compact.
There are a connected compact regular domain $\Gamma$,
$r\in C^\infty(M)$, and $c<0$ such that
\[
\begin{gathered}
 K=\{r\le0\},\qquad \Gamma=\{r\le c\},\\
 A\subset\intt\Gamma,\qquad \Gamma\subset\intt K.
\end{gathered}
\]
The pairs $(\Gamma,Y_\Gamma)$ and $(K,Y_K)$ are diffeomorphic,
and on $K\setminus\intt\Gamma$,
\begin{equation}\label{eq:collar-data}
 dr\ne0,\qquad
 \tr_{(\nabla_qr)^\perp}^q\Hess_qr>0.
\end{equation}
\end{lemma}

\begin{proof}
Choose a tubular neighborhood of $Y_K$ with signed
$q$-distance coordinate $s$, negative inside $K$.
Compactness and $H_q(Y_K)>0$ allow the collar to be chosen
so that its parallel levels are strictly mean convex: their second
fundamental forms depend smoothly on the collar parameter.
The smooth collar construction is standard
\cite[Theorem~9.25]{Lee}; strict mean convexity is preserved here by
compactness and continuity.
Shrink it so that it is disjoint from $A$.

Extend the signed coordinate to a smooth function $r$,
constant and negative in the deeper interior of $K$ and
constant and positive away from the exterior collar,
such that $K=\{r\le0\}$.
Choose $c<0$ sufficiently close to zero that
\[
 \Gamma=\{r\le c\}
\]
contains $A$ in its interior and $r=s$ throughout
$\{c\le r\le0\}$.
Then $\Gamma\subset\intt K$, and on the closed collar
$K\setminus\intt\Gamma$,
\[
 dr\ne0,\qquad
 \tr_{(\nabla_qr)^\perp}^q\Hess_qr>0.
\]
In particular, $Y_\Gamma$ is smooth and strictly mean convex.

An increasing change of the collar coordinate, equal to
the identity near its ends and sending $c$ to $0$, extends
by the identity to an ambient diffeomorphism taking
$(\Gamma,Y_\Gamma)$ to $(K,Y_K)$.
This proves connectedness and the asserted diffeomorphism.
\end{proof}

\begin{proof}[Proof of Lemma~\ref{lem:psc-collar}]
Apply Lemma~\ref{lem:collarcut} to the component $K$ in
Corollary~\ref{cor:psccomponent}.
The boundary diffeomorphism transfers a PSC metric from
$Y_K$ to $Y_\Gamma$.
Equation~\eqref{eq:collar-data} gives
\[
 H_q(Y_\Gamma)
 =\frac{\tr_{(\nabla_qr)^\perp}^q\Hess_qr}{|dr|_q}>0.
\]
Corollary~\ref{cor:boundaryepi} then gives
\[
 \pi_1(Y_\Gamma)\twoheadrightarrow\pi_1(\Gamma).
\]
\end{proof}

\subsection{Uniform protection for the relative maximum principle}

\begin{remark}[A quantitative rank-one correction]
Let $S$ be a symmetric form on a Riemannian 4-space, and let
$\xi\ne0$ be a covector. Suppose
\[
 \tr_{\ker\xi}^qS\ge\delta>0,\qquad
 \|S\|_{\mathrm{op},q}\le C,\qquad |\xi|_q\ge\alpha>0.
\]
Put $n=\xi^\sharp/|\xi|_q$.
For a unit normal $z=an+v$ to a three-plane $P$, with $v\perp n$
and $\tau=|v|_q^2=1-a^2$, we have
\[
\begin{aligned}
 \tr_P^q(S+L\xi\otimes\xi)
 &=\tr_{n^\perp}^qS+S(n,n)-S(z,z)+L|\xi|_q^2\tau\\
 &\ge\delta-2C\sqrt\tau-2C\tau+L\alpha^2\tau\\
 &\ge\frac\delta2+
       \left(L\alpha^2-2C-\frac{2C^2}{\delta}\right)\tau.
\end{aligned}
\]
The last step uses
$2C\sqrt\tau\le\delta/2+2C^2\tau/\delta$.
Consequently, for every three-plane $P$,
\begin{equation}\label{eq:rank-one-correction}
\begin{gathered}
 L\alpha^2\ge2C+\frac{2C^2}{\delta}\\
 \Longrightarrow\quad
 \tr_P^q(S+L\xi\otimes\xi)\ge\frac\delta2.
\end{gathered}
\end{equation}
For continuous data on a compact set with strict tangential trace
and nonvanishing $\xi$, the constants can be chosen uniformly.
\end{remark}

\begin{lemma}[Uniform protection from a mean-convex boundary]
\label{lem:bypass-protection}
Let $(M,q)$ be complete and connected, let $B$ be compact, and
suppose
\[
 \Gamma\subset\intt T,\qquad T\subseteq B,
\]
where $\Gamma$ is compact.
Given $\delta>0$, there is $a>0$, chosen independently of $D$ and $S$,
with the following property.
Let $D$ be a compact regular domain satisfying
\[
 T\subseteq D\subseteq B,\qquad H_q(Y_D)\ge\delta,
\]
and let $S\subseteq D$ be compact with
\[
 \mathsf{MP}_q(S;M\setminus\Gamma).
\]
Then
\[
 d_q(S,Y_D)\ge4a.
\]
\end{lemma}

\begin{proof}
Choose $r_0>0$ with
\[
 d_q(\Gamma,M\setminus\intt T)\ge r_0,
\]
and choose $\Lambda\ge0$ such that $\Ric_q\ge-\Lambda q$ on $B$.
Set
\[
 a=\frac14\min\left\{r_0,\frac{\delta}{\Lambda+1}\right\}>0.
\]
The assertion is immediate if $S=\varnothing$.

A point of $S\cap Y_D$ would lie outside $\Gamma$.
A local defining function for $D$, followed by a sufficiently
convex increasing reparametrization, would have positive Hessian
trace on every three-plane and a local maximum on $S$ there.
This follows from \eqref{eq:rank-one-correction} and
$H_q(Y_D)>0$, and contradicts the relative maximum principle.
Thus $S\subset\intt D$.

Realize
\[
 l=d_q(S,Y_D)>0
\]
at $x\in S$ and $y\in Y_D$.
If $x\in\Gamma$, then $l\ge r_0$.
Suppose henceforth that $x\notin\Gamma$.

Let $\gamma:[0,l]\to M$ be a minimizing geodesic from $y$ to $x$.
It starts in the inward normal direction and lies in $D$:
any later intersection with $Y_D$ would give a boundary point
closer to $x$.
We construct a smooth upper support for the boundary distance
without assuming that this distance is smooth at $x$.

Let $\Pi_s:T_yM\to T_{\gamma(s)}M$ be parallel transport.
For small $v\in T_xM$, write
\[
 v=\Pi_l(w+b\dot\gamma(0)),\qquad
 w\perp\dot\gamma(0),
\]
and define
\[
 a_0(w)
 =(\exp_y^q)^{-1}\bigl(\exp_y^{Y_D}(w)\bigr),
\]
where $\exp_y^{Y_D}$ is the intrinsic exponential map of $Y_D$.
Set
\[
 c_v(s)=
 \exp_{\gamma(s)}^q
 \left(
 \Pi_s\left[
 \left(1-\frac{s}{l}\right)a_0(w)
 +\frac{s}{l}\bigl(w+b\dot\gamma(0)\bigr)
 \right]
 \right).
\]
Then $c_v(0)\in Y_D$ and $c_v(l)=\exp_x^q(v)$.
For sufficiently small $v$, the curves have nonzero velocity,
so
\[
 \ell(\exp_x^q(v))=\operatorname{Length}_q(c_v)
\]
defines a smooth function near $x$.
It satisfies
\[
 d_q(\,\cdot\,,Y_D)\le\ell,\qquad
 \ell(x)=l,\qquad \nabla_q\ell(x)=\dot\gamma(l).
\]

For variations with $b=0$, the first-variation fields along
$\gamma$ are the parallel fields $\Pi_sw$.
The initial endpoint curves are intrinsic geodesics of $Y_D$,
and the final endpoint curves are ambient geodesics.
The first and second variation formulas
\cite[Lemma~5.4.2 and Theorem~6.1.4]{Petersen} therefore give
\[
 \tr_{\dot\gamma(l)^\perp}^q\Hess_q\ell(x)
 =
 -H_q(Y_D)(y)
 -\int_0^l\Ric_q(\dot\gamma,\dot\gamma)\,ds
 \le-\delta+\Lambda l.
\]
The boundary term is $-H_q$ because $\dot\gamma(0)$ is the
inward normal.
For these perpendicular parallel variation fields, the first
variation of speed vanishes pointwise, so the energy and length
second variations agree. The moving initial endpoint contributes
its second fundamental form; the final endpoint has zero
acceleration. No smoothness of the boundary distance is assumed.

On $S$ near $x$,
\[
 \ell\ge d_q(\,\cdot\,,Y_D)\ge l.
\]
Hence $e^{-L\ell}$ has a local maximum on $S$ at $x$.
If $\delta>\Lambda l$, then
\[
 \Hess_q(e^{-L\ell})
 =Le^{-L\ell}
   \bigl(-\Hess_q\ell+L\,d\ell\otimes d\ell\bigr)
\]
has positive trace on every three-plane for sufficiently large $L$.
Indeed, apply \eqref{eq:rank-one-correction} with
$S=-\Hess_q\ell$, $\xi=d\ell$, and tangential trace lower bound
$\delta-\Lambda l>0$.
This contradicts the relative maximum principle at $x$.

Thus, when $x\notin\Gamma$, one has $\delta\le\Lambda l$.
Together with the case $x\in\Gamma$, this gives
\[
 l\ge\min\left\{r_0,\frac{\delta}{\Lambda+1}\right\}=4a.
\]
\end{proof}

\subsection{Proof of Lemma~\ref{lem:bypass-capture}}

\begin{proof}
Assume $S\ne\varnothing$.
Put $\delta=\min_Bk>0$ and choose the protection radius $a$
from Lemma~\ref{lem:bypass-protection} before choosing any
inflation or forcing bump.
Since $H_q(Y_B)>k\ge\delta$, applying that lemma to $D=B$ gives
\[
 d_q(S,Y_B)\ge4a.
\]

Choose $0\le\beta\le1$ smooth, equal to one near $S$, with compact
support in
\[
 \{x:d_q(x,S)<a/2\}.
\]
Let $b_T$ be a smooth global defining function for $T$, with
\[
 T=\{b_T\le0\},\qquad Y_T=\{b_T=0\},
\]
and $db_T\ne0$ on $Y_T$.
Choose a regular value
\[
 t>\max\{0,\max_Sb_T\}
\]
of $b_T/\beta$ on $\{\beta>0\}$, and set
\[
 D=\{b_T-t\beta\le0\}.
\]
Sard's theorem applies to the smooth quotient on its open domain.
On its zero level with $\beta>0$,
\[
 d(b_T-t\beta)=\beta\,d(b_T/\beta)\ne0.
\]
At a zero with $\beta=0$, nonnegativity gives $d\beta=0$,
and hence
\[
 d(b_T-t\beta)=db_T\ne0.
\]
Thus $D$ is regular.
Moreover, $D\subset T\cup\supp\beta$, so it is compact, and
\[
 T\subseteq D,\qquad
 S\subset\intt D,\qquad
 D\subset\intt B.
\]
Outside $\supp\beta$, the domains $D$ and $T$ agree locally.

Choose $0\le\xi\le1$ smooth, equal to one near $\supp\beta$,
with compact support in $\{x:d_q(x,S)<a\}$.
Take
\[
 C>\max\{0,\max_{Y_D}(H_q(Y_D)-k)\},
\]
and set $\psi=k+C\xi$.
On the unchanged boundary, $H_q(Y_D)=k$; on the changed part,
$\xi=1$.
Therefore,
\[
 H_q(Y_D)\le\psi,\qquad
 \psi\ge k,\qquad
 \psi=k\text{ near }Y_B.
\]
The original barriers remain valid:
\[
 H_q(Y_A)<k\le\psi,\qquad H_q(Y_B)>k=\psi.
\]

Theorem~\ref{input:eichmair}, applied with this $\psi$, supplies
a compact regular domain $E$ such that
\[
 D\subseteq E\subset\intt B,\qquad H_q(Y_E)=\psi\ge\delta.
\]
In particular, $T\subseteq E$ and $S\subset\intt E$.
Applying Lemma~\ref{lem:bypass-protection} to $E$ with the same
radius $a$ gives
\[
 d_q(S,Y_E)\ge4a.
\]
Since $\supp\xi\subset\{x:d_q(x,S)<a\}$, we have $\xi=0$ on $Y_E$.
Thus
\[
 H_q(Y_E)=k.
\]
The weak outermost property of $T$ forces $E=T$, and consequently
$S\subset\intt T$.

The component $K$ is open and closed in the compact regular
domain $T$.
If $S\setminus K$ were nonempty, a maximum of $F$ on this compact
set would be a local maximum on $S$ outside $\Gamma$,
contradicting strict three-convexity.
Thus $S\subset K$.

If $S\not\subset\Gamma$, maximize $r$ on $S$.
At a maximizing point, $r>c$, and
\[
 \Hess_q(e^{Lr})
 =Le^{Lr}\bigl(\Hess_qr+L\,dr\otimes dr\bigr)
\]
has positive trace on every three-plane for sufficiently large $L$,
by \eqref{eq:rank-one-correction}.
This again contradicts the relative maximum principle.
Hence $S\subset\Gamma$.
\end{proof}

\subsection{The weak hull and a global defining level}

Choose a smooth convex function $\chi_a\ge|\cdot|$ equal to
$|\cdot|$ outside $(-a,a)$, and set
\[
 \operatorname{max}_a(s,t)=\frac{s+t+\chi_a(s-t)}2.
\]
Convexity gives $\chi_a(z)\le a$ for $|z|\le a$, and hence
\begin{equation}\label{eq:rounded-max-error}
\begin{gathered}
 0\le\operatorname{max}_a(s,t)-\max(s,t)\le\frac a2,\\
 |s-t|\ge a
 \quad\Longrightarrow\quad
 \operatorname{max}_a(s,t)=\max(s,t).
\end{gathered}
\end{equation}
For smooth functions $v,w$,
\begin{equation}\label{eq:bypass-rounded-max}
\begin{aligned}
 &\Hess_q\operatorname{max}_a(v,w)\\
 &\qquad=t\Hess_qv+(1-t)\Hess_qw
   +\frac{\chi_a''(v-w)}2\,d(v-w)\otimes d(v-w),
\end{aligned}
\end{equation}
where $t=(1+\chi_a'(v-w))/2\in[0,1]$.
Thus smooth maxima preserve weak three-convexity, and preserve
strict three-convexity when both inputs are strict.

\begin{lemma}[Maximum principle for the weak hull]
\label{lem:bypass-hull-MP}
For every nonempty compact $\Gamma$,
\[
 \mathsf{MP}_q(\widehat\Gamma^{\,q}_3;M\setminus\Gamma)
\]
holds.
If $F$ is a strictly three-convex exhaustion, the hull is compact.
\end{lemma}

\begin{proof}
The hull is closed and contains $\Gamma$.
If a strictly three-convex exhaustion $F$ exists, the hull is
contained in
\[
 \{F\le\max_\Gamma F\},
\]
and is therefore compact.

To prove the maximum principle, suppose that a smooth test
function $\psi$ has a local maximum on $\widehat\Gamma^{\,q}_3$
at
\[
 x\in\widehat\Gamma^{\,q}_3\setminus\Gamma
\]
and has positive Hessian trace on every three-plane at $x$.
Compactness of the Grassmannian and continuity of the Hessian
give a neighborhood with a uniform positive trace margin.
Subtracting a sufficiently small multiple of a local squared
distance, and then a constant, gives a smooth function $\psi_0$
on a relatively compact neighborhood
$\Omega\Subset M\setminus\Gamma$ such that
\[
\begin{gathered}
 \psi_0(x)=0,\\
 \lambda_{(3)}^q(\Hess_q\psi_0)\ge0
 \quad\text{on }\Omega,\\
 \psi_0\le-4d
 \quad\text{on }\widehat\Gamma^{\,q}_3\cap\partial\Omega
\end{gathered}
\]
for some $d>0$.

Every compact set outside the hull admits a smooth weakly
three-convex separator $u$ with $u\le-1$ on $\Gamma$ and $u>1$
on that compact set.
Indeed, choose finitely many separating functions
$u_1,\ldots,u_N$ such that
\[
 u_j\le-2\text{ on }\Gamma
\]
and the sets $\{u_j>2\}$ cover the given compact set.
Iterating smooth maxima with $(N-1)a/2<1$ gives the asserted
separator by \eqref{eq:rounded-max-error} and
\eqref{eq:bypass-rounded-max}.
For the empty set, use $u=-1$.

Apply this construction to the compact set
\[
 E=\{y\in\partial\Omega:\psi_0(y)\ge-2d\},
\]
which is disjoint from the hull.
Choose $A>0$ sufficiently large and put
\[
 v=\operatorname{max}_{d/4}(-d,Au-3d).
\]
Since $x$ belongs to the hull,
\[
 u(x)\le\max_\Gamma u\le-1.
\]
Thus $v=-d$ on $\Gamma$ and at $x$.
On $E$, choose $A$ large enough that
$\psi_0<v-d/4$.
On $\partial\Omega\setminus E$, this inequality follows from
$\psi_0<-2d$ and $v\ge-d$.
Hence
\[
 \psi_0<v-d/4\quad\text{on }\partial\Omega.
\]

Define a global function by taking
$\operatorname{max}_{d/4}(v,\psi_0)$ on $\Omega$ and $v$ outside.
The formulas agree near $\partial\Omega$, so the function is
smooth and weakly three-convex.
It equals $-d$ on $\Gamma$ and $0$ at $x$, contradicting the
definition of the hull.
\end{proof}

\begin{lemma}[Hull equality from compact capture]
\label{lem:bypass-hull-equality}
In Lemma~\ref{lem:bypass-capture}, suppose in addition that
\[
 F\le b\text{ on }\Gamma,\qquad F>b\text{ on }\partial B
\]
for some $b$.
Then $\widehat\Gamma^{\,q}_3=\Gamma$.
\end{lemma}

\begin{proof}
Put $S=\widehat\Gamma^{\,q}_3$.
Lemma~\ref{lem:bypass-hull-MP} gives its compactness and relative
maximum principle.
The defining inequality for $F$ gives $F\le b$ on $S$, so
$S\cap\partial B=\varnothing$.
Consequently,
\[
 S\setminus B=S\cap(M\setminus\intt B)
\]
is compact.
If this set were nonempty, a maximum of $F$ on it would be a
local maximum on $S$ outside $\Gamma$, contradicting the relative
maximum principle.
Thus $S\subset B$.
Lemma~\ref{lem:bypass-capture} now gives $S\subset\Gamma$.
The reverse inclusion follows from the definition.
\end{proof}

\begin{lemma}[A strict global defining level]
\label{lem:bypass-defining-level}
Let $\Gamma$ be a nonempty smooth compact regular domain with
$H_q(\partial\Gamma)>0$, and suppose
\[
 \widehat\Gamma^{\,q}_3=\Gamma.
\]
If $M$ has a smooth strictly three-convex exhaustion $F$, then
there is a smooth strictly three-convex exhaustion $G_0$ satisfying
\[
 \Gamma=\{G_0\le0\},\qquad
 \partial\Gamma=\{G_0=0\},
\]
with $dG_0\ne0$ on $\partial\Gamma$.
\end{lemma}

\begin{proof}
Choose a smooth global defining function $r$ with
\[
 \Gamma=\{r\le0\},\qquad
 \partial\Gamma=\{r=0\},
\]
with $dr\ne0$ on $\partial\Gamma$ and $r=1$ outside a compact set.
Positive mean curvature gives
\[
 \tr_{(\nabla_qr)^\perp}^q\Hess_qr
 =|dr|_qH_q(\partial\Gamma)>0.
\]
Apply \eqref{eq:rank-one-correction} uniformly on the compact
boundary.
By continuity, there are $L>0$ and $0<\epsilon<1$ such that
\[
 \lambda_{(3)}^q
 \bigl(\Hess_qr+L\,dr\otimes dr\bigr)>0
 \quad\text{on }U=\{|r|<\epsilon\}.
\]
Here $\epsilon$ is chosen small enough that $U$ lies in the
boundary neighborhood where this inequality holds.

Set
\[
 h=A(e^{Lr}-1),\qquad
 A=\frac4{1-e^{-L\epsilon}}.
\]
Then $h$ is strictly three-convex on $U$, has the sign of $r$, and
\[
 h|_{r=-\epsilon}=-4,\qquad h|_{r=0}=0,\qquad
 h|_{r=\epsilon}=4e^{L\epsilon}.
\]
Put
\[
 C=\{r\le\epsilon/2\},\qquad H=4e^{L\epsilon}+3.
\]

We construct a strictly three-convex exhaustion $v$ satisfying
\[
 -3<v\le-2\text{ on }\Gamma,\qquad
 v>H\text{ on }M\setminus\intt C.
\]
Choose a positive affine rescaling $f$ of $F$ with
\[
 f>-3\text{ on }M,\qquad f\le-2\text{ on }\Gamma.
\]
Choose a sufficiently large compact sublevel $P$ of $F$,
containing $C$, such that $f>H+4$ on $M\setminus P$.
The set
\[
 E=P\setminus\intt C
\]
is compact and disjoint from $\Gamma$.

Hull equality and the finite separator construction in
Lemma~\ref{lem:bypass-hull-MP} give a smooth weakly three-convex
function $w$ with
\[
 w\le-2\text{ on }\Gamma,\qquad w>2\text{ on }E.
\]
Adding a sufficiently small positive multiple of $F$ gives a
strictly three-convex function $u$ satisfying
\[
 u\le-1\text{ on }\Gamma,\qquad u>1\text{ on }E.
\]
The inequalities are preserved because $\Gamma\cup E$ is compact.

Choose $0<a_0<1$ and set
\[
 v=\operatorname{max}_{a_0}\bigl(f,(H+10)u-6\bigr).
\]
On $\Gamma$,
\[
 (H+10)u-6\le-H-16<f-1,
\]
so $v=f$ there.
Hence $-3<v\le-2$ on $\Gamma$.
On $E$, the second input is greater than $H+4$;
on $M\setminus P$, the first input is greater than $H+4$.
Thus $v>H$ on $M\setminus\intt C$.
Both inputs are strictly three-convex, and $v\ge f$, so $v$
is a strictly three-convex exhaustion.

For $a=1/10$, define
\[
 G_0=
 \begin{cases}
  \operatorname{max}_a(v,h),& |r|<\epsilon,\\
  v,& |r|\ge\epsilon.
 \end{cases}
\]
On $r=-\epsilon$, one has $v-h>1$.
On $r=\epsilon$, one has $v-h>3$.
Thus the two formulas agree on neighborhoods of the seams.

On $\partial\Gamma$, $h=0$ and $v\le-2$, so $G_0=h$
in a smaller collar.
For $-\epsilon<r<0$, if $h\le-1$, then
\[
 G_0\le\max(v,h)+\frac1{20}
 \le-1+\frac1{20}<0.
\]
If $-1<h<0$, then $h-v>1>a$, so $G_0=h<0$.
For $0<r<\epsilon$, one has $G_0\ge h>0$.
Outside $U$, the inner region has $v\le-2$, while the outer
region lies outside $C$, where $v>H>0$.
This proves
\[
 \Gamma=\{G_0\le0\},\qquad \partial\Gamma=\{G_0=0\}.
\]
Since $G_0=h$ near the boundary, $dG_0\ne0$ there.

Equation~\eqref{eq:bypass-rounded-max} gives strict three-convexity,
and $G_0\ge v$ gives boundedness below and compact sublevels.
\end{proof}

\begin{lemma}[Morse approximation retaining the defining level]
\label{lem:bypass-relative-Morse}
The function $G_0$ in Lemma~\ref{lem:bypass-defining-level}
can be replaced by a smooth strictly three-convex Morse
exhaustion $G$ with the same zero sublevel and the same germ
near $\partial\Gamma$.
Every critical point of $G$ has index at most two.
\end{lemma}

\begin{proof}
Choose $0\le\beta\le1$ smooth, equal to $1$ near $\{G_0=0\}$,
with compact support contained in $\{dG_0\ne0\}$.
On that support,
\[
 |dG_0|_q^2\ge m_0>0.
\]
Properness gives a positive lower bound for $|G_0|$ outside
$\intt\{\beta=1\}$.

Apply \cite[Chapter~6, Section~1, Theorem~1.2]{Hirsch} with
$s=\infty$. The prescribed pointwise tolerances on the value,
differential, and Hessian define a strong $C^\infty$ neighborhood,
so this gives smooth Morse approximants satisfying the required
fine $C^2$ bounds.
Choose a smooth Morse function $H$ with sufficiently fine
value, differential, and Hessian approximation to $G_0$.
Put $E=H-G_0$ and
\[
 G=G_0+(1-\beta)E.
\]
Then $G=G_0$ near the zero level and $G=H$ outside $\supp\beta$.
The identity
\[
 dG(\nabla_qG_0)
 =|dG_0|_q^2
  +(1-\beta)dE(\nabla_qG_0)
  -E\,d\beta(\nabla_qG_0)
\]
keeps this expression positive on $\supp\beta$ when the first
jets are sufficiently close.
Thus every critical point lies where $G$ agrees locally with
the Morse function $H$.

Choose a continuous function $\tau>0$ such that
\[
 \tr_P^q\Hess_qG_0\ge\tau(x)
\]
for every three-plane at $x$.
Let $A_1,A_2$ bound $|d\beta|_q$ and
$\|\Hess_q\beta\|_{\mathrm{op},q}$.
If
\[
 \max\bigl\{|E|,\ |dE|_q,\ \|\Hess_qE\|_{\mathrm{op},q}\bigr\}
 \le\eta(x),
\]
then
\begin{align*}
 \Hess_q(G-G_0)
 &=(1-\beta)\Hess_qE-E\Hess_q\beta\\
 &\quad-d\beta\otimes dE-dE\otimes d\beta,\\
 \left|\tr_P^q\Hess_q(G-G_0)\right|
 &\le3\eta(x)(1+2A_1+A_2).
\end{align*}
Choose the fine approximation so that this last quantity is
less than $\tau(x)$.
This preserves strict three-convexity.

A sufficiently small value error preserves both signs outside
$\intt\{\beta=1\}$, while $G=G_0$ on $\{\beta=1\}$.
Taking also $|G-G_0|<1$ preserves boundedness below and compact
sublevels.
Thus $G$ has the same zero sublevel and the same germ near
$\partial\Gamma$.

Finally, three negative critical Hessian eigenvalues would
contradict the positive three-plane trace, so
\[
 \ind_xG\le2\qquad(x\in\Crit G).
\]
\end{proof}

\subsection{Proof of Theorem~\ref{thm:psc-defining}}

\begin{proof}
Retain the fixed potential, metric $q$, exhaustion $F$, and
regular cores from
Theorems~\ref{thm:joint}--\ref{thm:cores}.
Given a compact set $Q\subset M$, choose $i$ such that
\[
 Q\subset\intt C_i.
\]
Set $A=C_i$, and take $B=C_j$ and the forcing function $k$
from Proposition~\ref{prop:forcing}.

Proposition~\ref{prop:outerpsc} gives the weakly outermost
domain $T$, and Corollary~\ref{cor:psccomponent} gives its
component $K$ containing $A$.
Its boundary is connected and admits PSC, and inclusion
surjects on fundamental groups.

Lemmas~\ref{lem:psc-collar} and~\ref{lem:collarcut} give a
connected inward collar cut $\Gamma$ with
\[
 Q\subset\intt A\subset\intt\Gamma,\qquad
 \Gamma\subset\intt K,
\]
positive outward mean curvature, and collar data $r,c$
as required in Lemma~\ref{lem:bypass-capture}.
The boundary is diffeomorphic to $Y_K$, hence admits PSC.
Corollary~\ref{cor:boundaryepi} gives connectedness and
\[
 \pi_1(Y_\Gamma)\twoheadrightarrow\pi_1(\Gamma).
\]

Since $B=C_j$ is a regular component of $\{F\le s_j\}$,
every point of $B$ with $F=s_j$ is a boundary point.
Thus $F<s_j$ on $\Gamma$.
Compactness gives a number $b$ such that
\[
 \max_\Gamma F<b<s_j=F|_{\partial B}.
\]
Lemmas~\ref{lem:bypass-capture} and
\ref{lem:bypass-hull-equality} imply
\[
 \widehat\Gamma^{\,q}_3=\Gamma.
\]
This also proves Lemma~\ref{lem:bypass-hull-body}.

Finally, Lemmas~\ref{lem:bypass-defining-level} and
\ref{lem:bypass-relative-Morse} give a smooth strictly
three-convex Morse exhaustion $G$ with
\[
 \Gamma=\{G\le0\},\qquad Y_\Gamma=\{G=0\},
\]
and $dG\ne0$ on $Y_\Gamma$.
Its critical indices are at most two.
The construction keeps $q$ fixed while $Q$ varies and
establishes all the asserted properties.
\end{proof}

\Needspace{10\baselineskip}
\section{Proofs of Section~\ref{sec:topology}}
\label{app:topology}

All homology groups in this section have integer coefficients
unless otherwise specified.

\subsection{Proof of Lemma~\ref{lem:bandorientation}}

\begin{proof}
Let $a$ realize $\mathsf{Band}_q(K,B)$, and let $D$ denote
either $K$ or $B$.

Suppose $x\in Y_D$ and $da_x=0$.
Since $a|_{Y_D}$ is constant,
\[
 (\Hess_q a)_x(V,W)=0
 \qquad(V,W\in T_xY_D).
\]
In a frame adapted to $T_xY_D$, the Hessian has the form
\[
 \begin{pmatrix}
  0_{3\times3}&\beta\\
  \beta^{\mathsf T}&\gamma
 \end{pmatrix},
\]
whose rank is at most two.
This contradicts the nondegeneracy in
Definition~\ref{def:band}.
Hence $da\ne0$ on both boundaries.

Let $\nu_q$ be the outward unit normal of $D$.
At a boundary point,
\[
 \nabla_q a=(\nu_qa)\nu_q.
\]
For tangent vectors $V,W\in TY_D$, we also have
\[
 \Hess_q a(V,W)=(\nu_qa)\mathrm{II}_q(V,W).
\]
Therefore
\[
 0<
 \tr_{(\nabla_q a)^\perp}^q\Hess_q a
 =(\nu_qa)H_q(Y_D).
\]
Since $H_q(Y_D)>0$, we have $\nu_qa>0$.
Thus $\nabla_q a$ points outward from both $K$ and $B$.
\end{proof}

\subsection{Proof of Proposition~\ref{prop:transfer}}
\label{app:transfer}

\begin{lemma}[Relative handles for a band and its tail]
\label{input:bandtail}
Suppose $K\subset\intt B$ are compact regular domains,
$\mathsf{Band}_q(K,B)$ holds, and its boundary gradients
point outward from $K$ and $B$.
Let $F$ be a proper, bounded-below Morse function whose
critical indices are at most two.
Suppose $B$ is a component of a regular sublevel of $F$.

Then every finite CW structure on $Y_K$ extends to a
relative cell model of $(E_K,Y_K)$ with relative cells
only in dimensions $0,1,2$.
Moreover,
\begin{equation}\label{eq:relativehomology}
\begin{gathered}
 H_\ell(M,K)=0\quad(\ell\ge3),\\
 H_2(M,K)\text{ is torsion-free}.
\end{gathered}
\end{equation}
\end{lemma}

\begin{proof}
Let $a$ realize $\mathsf{Band}_q(K,B)$, and put
\[
 W=B\setminus\intt K.
\]
We first replace $a|_W$ by a Morse function with range $[0,1]$
and the same critical points and indices.

Choose $0<\varepsilon<1/2$ sufficiently small that
\[
 \frac14<\varepsilon a+\frac12<\frac34
 \quad\text{on }W.
\]
The outward boundary gradients make $a$ a normal coordinate
in disjoint boundary collars.
On the $W$-side of $Y_K$, we have $a>0$ near the boundary;
on the $W$-side of $Y_B$, we have $a<1$.

Shrink the collars so that their $a$-coordinate ranges are
$[0,\delta]$ and $[1-\delta,1]$, where $0<\delta<1/2$.
Choose nonnegative smooth bumps $\eta_-,\eta_+$ supported
in the interiors of these intervals, each with integral one.
Define
\[
\begin{aligned}
 r_-(s)&=\varepsilon s+
          \frac12\int_0^s\eta_-(t)\,dt,\\
 r_+(s)&=\varepsilon s+\frac12+
          \left(\frac12-\varepsilon\right)
          \int_{1-\delta}^s\eta_+(t)\,dt.
\end{aligned}
\]
Then $r_-(0)=0$, $r_+(1)=1$, and
\[
\begin{aligned}
 r_-'&=\varepsilon+\tfrac12\eta_->0,\\
 r_+'&=\varepsilon+
        (\tfrac12-\varepsilon)\eta_+>0.
\end{aligned}
\]
Near the interior edge of either collar, the corresponding
function agrees with $\varepsilon s+1/2$.
Use $r_-\circ a$ and $r_+\circ a$ in the two collars and
$\varepsilon a+1/2$ elsewhere.
The resulting smooth function
\[
 \widehat a:W\longrightarrow[0,1]
\]
satisfies
\[
 \widehat a^{-1}(0)=Y_K,\qquad
 \widehat a^{-1}(1)=Y_B.
\]
It introduces no critical points in the collars.
Elsewhere it is a positive affine transformation of $a$,
so its critical points and indices agree with those of $a$.

The relative Morse construction builds $W$ from a collar
of $Y_K$ by handles of indices $0,1,2$.
Each handle contributes a relative cell of its index
\cite[Theorems~3.13--3.14]{MilnorHC}.

Let $s$ be the regular value defining $B$, and write
\[
 \{F\le s\}=B\sqcup B_1\sqcup\cdots\sqcup B_\ell.
\]
There are finitely many components because the sublevel
is a compact smooth manifold with boundary.
Each $B_\alpha$ is built from the empty set by the
$F$-handles in that component, again of indices $0,1,2$.

Choose regular values above $s$ tending to infinity.
Each intervening slab is compact and contains finitely
many critical points, so only finitely many handles are
attached at each stage.
Critical points sharing a critical value can be treated
simultaneously
\cite[remark following Theorem~3.14]{MilnorHC}.
Together with the compact band and the components
$B_\alpha$, these handles construct $E_K$ relative to $Y_K$.

Start with the prescribed finite CW structure on $Y_K$.
Replace handles by core disks and apply cellular approximation
to their attaching maps \cite[Theorem~4.8]{Hatcher}.
Homotopic attaching maps give homotopy equivalent adjunction
spaces relative to the preceding stage
\cite[Proposition~0.18]{Hatcher}.
Choose the finite-stage models and comparison maps inductively
relative to the preceding stage and to $Y_K$.
Their union gives a CW complex $X$, containing the prescribed
$Y_K$ as a subcomplex, and a comparison map
\[
 e:X\longrightarrow E_K
\]
whose restriction to $Y_K$ is homotopic to its inclusion.

Every compact subset of $E_K$ lies in a finite stage, and
every compact subset of $X$ lies in a finite subcomplex
\cite[Proposition~A.1]{Hatcher}.
Consequently, maps of spheres and their homotopies are detected
at finite stages.
The finite-stage equivalences imply that $e$ is bijective on
components and induces isomorphisms on all homotopy groups.
Since smooth manifolds have CW type, Whitehead's theorem,
applied componentwise, makes $e$ a homotopy equivalence
\cite[Theorem~4.5]{Hatcher}.

The homotopy extension property of the CW pair $(X,Y_K)$
allows $e$ to be replaced by a homotopic map satisfying
\[
 e|_{Y_K}=(Y_K\hookrightarrow E_K)
\]
exactly \cite[Proposition~0.16]{Hatcher}.
This gives the relative model of Definition~\ref{def:CW},
including its prescribed boundary cells and actual boundary map.

Collared excision gives
\[
 H_\ell(M,K)\cong H_\ell(E_K,Y_K).
\]
Indeed, excising a slightly smaller core of $K$ leaves a collar
pair which deformation retracts to $(E_K,Y_K)$.
Since $e$ is a homotopy equivalence and restricts to the actual
inclusion of $Y_K$, the pair sequences give
\[
 H_\ell(X,Y_K)\cong H_\ell(E_K,Y_K).
\]
Relative cellular homology
\cite[Theorem~2.35 and Section~2.2]{Hatcher} computes these
groups using a chain complex with
\[
 C_\ell(X,Y_K)=0\quad(\ell\ge3).
\]
Hence
\[
\begin{gathered}
 H_\ell(M,K)=0\quad(\ell\ge3),\\
 H_2(M,K)\cong
 \ker\bigl(\partial_2:C_2(X,Y_K)\to C_1(X,Y_K)\bigr).
\end{gathered}
\]
The latter is a subgroup of a free abelian group and is
therefore torsion-free.
\end{proof}

\begin{proof}[Proof of Proposition~\ref{prop:transfer}]
Put $B=D_j$ and $W=B\setminus\intt K$.
The boundary $Y_B$ is strictly mean convex by
Theorem~\ref{thm:cores}.
Together with $H_q(Y_K)>0$, Lemma~\ref{lem:bandorientation}
supplies the outward boundary gradients.
Apply Lemma~\ref{input:bandtail} with tail function $F=H$.
It gives the relative cell model and the homology assertions.

To identify the exterior fundamental-group map, read the
compact band from $Y_B$ toward $Y_K$.
A forward handle of index $\lambda\le2$ becomes a reverse
handle of index $4-\lambda\ge2$.
Thus $W$ is built from a collar of the connected $Y_B$
without relative zero- or one-handles.
It follows that $W$ is path connected and
\[
 \pi_1(Y_B)\twoheadrightarrow\pi_1(W).
\]

The closed exteriors satisfy
\[
 E_K=E_B\cup W,\qquad E_B\cap W=Y_B.
\]
Boundary collars give an open cover with these homotopy types.
Moreover, $E_B\simeq O_B$ is path connected by
Theorem~\ref{thm:cores}.
Van Kampen \cite[Theorem~1.20]{Hatcher} gives
\[
 \pi_1(E_K)
 \cong\pi_1(E_B)*_{\pi_1(Y_B)}\pi_1(W).
\]
Since $\pi_1(Y_B)\to\pi_1(W)$ is surjective, the first
vertex group surjects onto this pushout:
\[
 \pi_1(E_B)\twoheadrightarrow\pi_1(E_K).
\]

The commuting square of actual inclusions
\[
\begin{array}{ccc}
 O_B&\longrightarrow&O_K\\
 \big\downarrow&&\big\downarrow\\
 E_B&\longrightarrow&E_K
\end{array}
\]
has homotopy equivalences as its vertical maps.
With compatible basepoints, it follows that
\[
 \pi_1(O_{D_j})\twoheadrightarrow\pi_1(O_K)
\]
is induced by inclusion.
The same decomposition shows that $E_K$ is path connected,
and hence so is $O_K$.
\end{proof}

\subsection{Proof of Lemma~\ref{lem:boundaryH}}

\begin{proof}
Contractibility of $M$ implies orientability of $M$ and $K$.
Since $K$ is compact, its homology groups are finitely generated.
Thus the torsion-free group $H_1(K)$ is free abelian.

The universal coefficient theorem
\cite[Theorem~3.2]{Hatcher} gives
\[
 H^2(K)=H^3(K)=0.
\]
Indeed, the relevant terms are
\[
\begin{gathered}
 \Hom(H_2(K),\Z),\quad \Ext(H_1(K),\Z),\\
 \Hom(H_3(K),\Z),\quad \Ext(H_2(K),\Z),
\end{gathered}
\]
all of which vanish.
Poincar\'e--Lefschetz duality
\cite[Theorem~3.43]{Hatcher} therefore gives
\[
\begin{aligned}
 H_2(K,Y_K)&\cong H^2(K)=0,\\
 H_1(K,Y_K)&\cong H^3(K)=0.
\end{aligned}
\]
The exact sequence of the pair shows that inclusion induces
\[
 H_1(Y_K)\xrightarrow{\cong}H_1(K).
\]

Use boundary collars to apply Mayer--Vietoris to
$M=K\cup E_K$.
Since $M$ is contractible, the inclusion maps give isomorphisms
\[
 H_2(Y_K)\xrightarrow{\cong}
 H_2(K)\oplus H_2(E_K)
\]
and
\[
 H_1(Y_K)\xrightarrow{\cong}
 H_1(K)\oplus H_1(E_K),
\]
up to the usual choice of signs.
As $H_2(K)=0$, the first proves the asserted isomorphism
$H_2(Y_K)\to H_2(E_K)$.

In the second, projection to $H_1(K)$ is already an isomorphism.
Surjectivity onto an element $(0,\xi)$ therefore forces $\xi=0$.
Hence $H_1(E_K)=0$.
Finally, $H_1(Y_K)\cong H_1(K)$ is torsion-free.
\end{proof}

\subsection{Proof of Proposition~\ref{input:normalform}}

\begin{proof}
The simply connected ambient manifold is orientable, so
$K$ and $Y=\partial K$ inherit orientations.
The PSC classification
\cite[Proposition~2.2]{CMM}, applied to the identity map of $Y$,
gives
\[
 Y\cong
 \left(\mathop{\#}_{\ell=1}^{v}P_\ell\right)
 \#\left(\mathop{\#}_{i=1}^{r}(S^1\times S^2)\right),
\]
where the $P_\ell$ are spherical three-manifolds.
An empty connected sum is interpreted as $S^3$.

Connected-sum homology gives
\[
 H_1(Y)\cong
 \left(\bigoplus_{\ell=1}^{v}H_1(P_\ell)\right)\oplus\Z^r.
\]
Each $H_1(P_\ell)$ is finite.
Since $H_1(Y)$ is torsion-free, all these finite summands vanish.
Thus each $P_\ell$ is an integral homology sphere.
By \cite[Corollary~2.3]{CMM}, it is $S^3$ or a connected sum
of copies of $P$ and $-P$.
Its fundamental group is finite, whereas a free product of two
nontrivial groups is infinite.
It therefore has at most one nontrivial summand.
After discarding the $S^3$ summands, each $P_\ell$ is $P$ or $-P$.
Put
\[
 Y_0=\mathop{\#}_{\ell=1}^{v}P_\ell.
\]

For each punctured $P_\ell$, choose a finite PL spine of
dimension at most two.
Theorem~2.5 of \cite{Spines} supplies a special spine of this
compact 3-manifold with boundary, and Theorem~2.1 identifies
the manifold with a regular neighborhood of its spine and
its spine complement with a boundary collar.
For each $S^1\times S^2$ summand, use the product CW structure
with embedded two-skeleton
\[
 (S^1\times\{p_i\})\cup(\{t_i\}\times S^2)
 \cong S^1\vee S^2.
\]
Its complement is an open three-cell:
\[
 (S^1\setminus\{t_i\})\times(S^2\setminus\{p_i\})
 \cong\R^3.
\]
The circle contributes one one-cell, and the sphere contributes
one two-cell with constant attaching map at the wedge point.

Join the punctured $P_\ell$ spines and the marked wedge points
by a tree of arcs, meeting each spine at a single chosen point.
Let $A$ be the union of the $P_\ell$ spines and this tree.
Choose finite CW structures on the spines and tree compatible
at these attaching points.
The punctured homology spheres are acyclic, and the connecting
tree adds no reduced homology. Thus
\[
 \widetilde H_*(A)=0,\qquad
 \pi_1(A)\cong\mathop{*}_{\ell=1}^{v}\pi_1(P_\ell).
\]
If $v=0$, take $A$ to be the connecting tree, or a point if
no connecting arcs are needed.
Every selected sphere cell is attached at a point of $A$.

Finite groups are residually finite.
Gruenberg's theorem that free products of residually finite
groups are residually finite, recalled in
\cite[Section~1.1]{Berlai}, therefore makes $\pi_1(A)$
residually finite.

Take regular neighborhoods of the summand spines and join them
by boundary bands along the chosen tree.
Their boundary connected sum is the punctured connected sum
$Y\setminus\intt B^3$, with the joined spine as a spine.
The regular-neighborhood description identifies its complement
with a boundary collar.
Capping by $B^3$ makes the complement of the joined spine an
open three-ball. The regular-neighborhood collapse on its
boundary supplies a characteristic map for this three-cell,
and hence extends the spine structure to a CW structure on $Y$
with one three-cell.
Every two-cell outside $A$ is one of the selected sphere cells;
there is one for each $S^1\times S^2$ summand, so $s=r$.
The circle one-cells supply the free generators outside
$\pi_1(A)$.

Now choose disjoint loops
\[
 \ell_i=S^1\times\{p_i\}\qquad(1\le i\le r)
\]
in the $S^1\times S^2$ summands, away from the connected-sum
balls.
Attach 4-dimensional two-handles to $K$ along these loops
with their product framings.
After rounding corners, denote the resulting compact smooth
4-manifold by $D$.

On the boundary, each surgery replaces $S^1\times D^2$
by $D^2\times S^1$.
The complement of the removed solid torus in its
$S^1\times S^2$ summand is another solid torus, and the
product-framed gluing gives
\[
 (S^1\times D^2)\cup(D^2\times S^1)
 \cong\partial(D^2\times D^2)\cong S^3.
\]
Thus each surgery replaces one $S^1\times S^2$ summand by $S^3$,
and
\[
 \partial D\cong Y_0.
\]
This is the standard product-framed surgery on the $S^1$ factor
of $S^1\times S^2$, modeled by the index-two handle cobordism
$W_2$ in \cite[Example~8.26(i)]{Ranicki}, with $m=3$ and $i=1$.
In that model, the preceding one-handle and the two-handle cancel
by \cite[Lemma~8.24]{Ranicki}.
Here the two-handles are attached abstractly to $K$, so no
spanning disks in the original ambient manifold are required.

The handle attachment and corner rounding give a continuous
embedding $f:K\hookrightarrow D$.
Placing $D$ in its smooth double realizes it as a compact
regular domain in an auxiliary smooth 4-manifold.
Its boundary is connected, and
\[
\begin{aligned}
 \pi_1(\partial D)&\cong\pi_1(Y_0)\cong\pi_1(A),\\
 \pi_1(Y)&\cong\pi_1(\partial D)*F_r.
\end{aligned}
\]
Moreover, $H_1(Y)\cong\Z^r$ has basis
\[
 b_i=[\ell_i]\qquad(1\le i\le r).
\]

The pair $(D,K)$ is obtained by attaching exactly $r$ two-handles.
Replacing the handles by their core disks gives
\[
 H_1(D,K)=0,\qquad H_2(D,K)\cong\Z^r.
\]
Let $e_i$ be the relative class of the oriented core disk of
the $i$th handle.
Its boundary is the attaching loop $\ell_i$, so the actual
connecting homomorphism satisfies
\[
 \delta e_i=(Y\hookrightarrow K)_*[\ell_i]
           =(Y\hookrightarrow K)_*b_i.
\]
Choose the core-disk orientations consistently with the loops.
All the cell, group, and relative-homology identifications
therefore arise from the same construction.
\end{proof}

\subsection{Proof of Theorem~\ref{thm:unitary}}

\begin{proof}
Put $H=\pi_1(A)$.
Choose a maximal tree in the one-skeleton of $Z/A$.
Its edges lift to a forest in $Z$ whose components each meet
$A$ in exactly one vertex.
Contract each lifted tree to that vertex.
This fixes $A$ and preserves the homotopy type of the pair,
so we may assume there are no relative zero-cells.

Let $x_1,\ldots,x_n$ correspond to the remaining relative
one-cells, and let $w_1,\ldots,w_m$ be the attaching words
of the relative two-cells.
The relative cellular presentation is
\[
 \pi_1(Z)\cong
 \bigl(H*F(x_1,\ldots,x_n)\bigr)/
 \langle\!\langle w_1,\ldots,w_m\rangle\!\rangle.
\]
Each $w_j$ has coefficients in $H$, and the map from $H$
is the original inclusion-induced map
\cite[Proposition~1.26(a)]{Hatcher}.

The relative cellular chain complex is
\[
 0\longrightarrow\Z^m
 \xrightarrow{\mathcal D}\Z^n
 \longrightarrow0.
\]
The entry $\mathcal D_{ij}$ is the exponent sum of $x_i$
in $w_j$ \cite[Theorem~2.35 and Section~2.2]{Hatcher}.
Since $H_2(Z,A)=0$, the matrix $\mathcal D$ is injective
over $\Z$ and has rank $m$ over $\Q$.

If $m=0$, assign arbitrary unitary matrices to the new generators.
Otherwise, choose $m$ variables corresponding to a nonzero
$m\times m$ minor of $\mathcal D$ and assign the identity
to all other variables.
Apply $\rho$ to the coefficients in the attaching words.
This gives $m$ equations in $m$ variables in $U(d)$ whose
exponent-sum determinant is nonzero.

Since $U(d)$ is a compact connected Lie group,
\cite[Theorem~1(ii)]{GR} gives a solution in this same $U(d)$.
Equivalently, one may apply \cite[Lemma~2.1(1)]{NitscheThom}
directly to the original $n$-variable system:
after erasing the coefficients, its presentation complex has
cellular two-boundary $\mathcal D$, and hence vanishing $H_2$.
The resulting homomorphism on $H*F(x_1,\ldots,x_n)$
annihilates every attaching word and restricts to $\rho$ on $H$.
It therefore descends to
\[
 \widetilde\rho:\pi_1(Z)\longrightarrow U(d),
 \qquad \widetilde\rho\circ i_*=\rho.
\]

Suppose now that $H$ is residually finite, and let
$1\ne\gamma\in H$.
Choose a homomorphism $b:H\to Q$ to a finite group with
$b(\gamma)\ne1$.
Its regular permutation representation is faithful and unitary:
\[
 \lambda:Q\hookrightarrow U(|Q|).
\]
Extend $\lambda\circ b$ by the preceding argument.
If $i_*(\gamma)=1$, the extension would imply
$\lambda(b(\gamma))=I$, a contradiction.
Hence $i_*$ is injective.
\end{proof}

\subsection{Proof of Corollary~\ref{cor:spine}}

\begin{proof}
Suppose a loop in $A$ is nullhomotopic in $Z$.
Its nullhomotopy has compact image and is therefore contained,
together with the finite subcomplex $A$, in a finite subcomplex
$Z_0\subset Z$ \cite[Proposition~A.1]{Hatcher}.
Taking the component containing $A$, we may assume $Z_0$
is connected.

Since $Z$ has dimension at most two, cellular chains give
\[
 H_2(Z_0)\hookrightarrow H_2(Z)=0.
\]
The exact sequence
\[
 H_2(Z_0)\longrightarrow H_2(Z_0,A)\longrightarrow H_1(A)
\]
and $H_1(A)=0$ imply $H_2(Z_0,A)=0$.
Theorem~\ref{thm:unitary} makes
$\pi_1(A)\to\pi_1(Z_0)$ injective.
Thus the loop is already nullhomotopic in $A$, proving
injectivity into $\pi_1(Z)$.
The finiteness assertion follows immediately.
\end{proof}

\subsection{Proof of Lemma~\ref{lem:deleted}}

\begin{proof}
The selected two-cells have constant attaching maps at points
of $A$.
Deleting their open cells from $X^{(2)}$ therefore leaves a
subcomplex $Z$ containing $j(A)$.
It retains the entire one-skeleton of the connected complex
$X$, so $Z$ is connected.

Write $\partial_i^T$ for the cellular boundary maps of a
CW complex $T$.
The closed connected orientable three-manifold $Y$ has one
three-cell and no higher cells.
Since $H_3(Y)=\Z$, its cellular differential satisfies
\[
 \partial_3^Y=0,\qquad
 H_2(Y^{(2)})\xrightarrow{\cong}H_2(Y).
\]
The relative model adds no cells of dimension at least three.
Thus
\[
 C_3(X)=j_*C_3(Y),\qquad \partial_3^X=0,
\]
and
\[
 H_2(X^{(2)})\xrightarrow{\cong}H_2(X).
\]

Every two-cell of $Y$ other than the selected sphere cells
belongs to $A$.
Since $\dim A\le2$ and $H_2(A)=0$, the differential
$\partial_2^A$ is injective.
The selected cells have zero cellular boundary, so
\[
 \ker\partial_2^Y
 =\bigoplus_{\alpha=1}^{s}\Z[e^2_\alpha].
\]
The additional one-cells outside $A$ do not affect this
calculation: boundaries of two-cells in $A$ lie in $C_1(A)$.

By the defining identity
\[
 e\circ j=(Y\hookrightarrow E_K)
\]
and Lemma~\ref{lem:boundaryH}, the map
\[
 H_2(Y^{(2)})
 \xrightarrow{\,j_*\,}H_2(X^{(2)})
\]
is an isomorphism.
Indeed, both skeleton inclusions are isomorphisms in degree
two, $e$ is a homotopy equivalence, and the actual inclusion
$Y\hookrightarrow E_K$ induces an isomorphism on $H_2$.
As the cellular three-boundaries vanish, this identifies
the cycle groups themselves:
\[
 \ker\partial_2^X
 =\bigoplus_{\alpha=1}^{s}\Z[j(e^2_\alpha)].
\]

On the other hand,
\[
 C_2(X)=C_2(Z)\oplus
        \bigoplus_{\alpha=1}^{s}\Z[j(e^2_\alpha)].
\]
Each selected cell has zero boundary, and hence
\[
 \ker\partial_2^X
 =\ker\partial_2^Z\oplus
  \bigoplus_{\alpha=1}^{s}\Z[j(e^2_\alpha)].
\]
Comparison gives $\ker\partial_2^Z=0$.
Since $Z$ has dimension at most two, $H_2(Z)=0$.

Adding back two-cells with constant attaching maps does not
change the fundamental group.
Neither does adding the boundary three-cell.
Consequently,
\[
 \pi_1(Z)\cong\pi_1(X^{(2)})
 \cong\pi_1(X)\cong\pi_1(E_K).
\]
The equivalence $E_K\simeq O_K$ identifies this group with
$\pi_1(O_K)$.

The subcomplex $j(A)$ is finite and connected, with
$H_1(j(A))=0$ and residually finite fundamental group.
Corollary~\ref{cor:spine} gives
\[
 \pi_1(A)\hookrightarrow\pi_1(Z)
 \cong\pi_1(E_K)\cong\pi_1(O_K).
\]
This is the asserted injection.
\end{proof}

\subsection{Proof of Proposition~\ref{prop:freeboundary}}

\begin{proof}
Lemma~\ref{lem:deleted} gives an injection
\[
 \pi_1(A)\hookrightarrow\pi_1(O_K).
\]
By Proposition~\ref{input:normalform},
\[
 \pi_1(A)\cong\mathop{*}_{\ell=1}^{v}I^*,
 \qquad I^*=\pi_1(P),
\]
where $I^*$ is the binary icosahedral group.

We first show that $v\le1$ using the group hypothesis in
Condition~\ref{cond:topology}. If $v\ge2$, choose three distinct
elements $a_0=1,a_1,a_2$ in the first factor and an involution $b$
in the second factor. Such an involution is the central element
of order two in $I^*$.
The free-product normal form shows that
\[
 t_i=a_i b a_i^{-1}\quad(0\le i\le2)
\]
generate a subgroup isomorphic to
$(\Z/2)*(\Z/2)*(\Z/2)$: a word in the $t_i$ with unequal adjacent
indices expands to a nonempty reduced word in the two original
factors. The kernel of the homomorphism sending all three $t_i$
to $1\in\Z/2$ is free of rank two. Indeed, Schreier rewriting
with coset representatives $1,t_0$ gives the free basis
$t_0t_1,t_0t_2$ and no nontrivial relators.
Thus $F_2$ embeds in $\pi_1(A)$ and hence in $\pi_1(O_K)$,
contrary to Condition~\ref{cond:topology}. Therefore $v\le1$.

Lemma~\ref{lem:boundaryH} gives an inclusion-induced isomorphism
\[
 i_*:H_1(Y)\xrightarrow{\cong}H_1(K).
\]
For the bases in \eqref{eq:actualdelta}, we have
$\delta e_i=i_*b_i$.
Therefore the actual connecting homomorphism is an isomorphism:
\[
 \delta:H_2(D,K)\xrightarrow{\cong}H_1(K).
\]
The exact sequence
\[
\begin{aligned}
 0=H_2(K)&\longrightarrow H_2(D)
 \longrightarrow H_2(D,K)
 \xrightarrow{\delta}H_1(K)\\
 &\longrightarrow H_1(D)
 \longrightarrow H_1(D,K)=0
\end{aligned}
\]
now gives
\[
 H_1(D)=H_2(D)=0.
\]

The universal coefficient theorem
\cite[Theorem~3.2]{Hatcher} implies
\[
 H^1(D;\Z/2)=H^2(D;\Z/2)=0.
\]
Hence the obstruction classes
\[
 w_1(TD)\in H^1(D;\Z/2),\qquad
 w_2(TD)\in H^2(D;\Z/2)
\]
vanish.
Thus $D$ is orientable and spin.
Moreover, $H_2(D;\R)=0$, so $\sigma(D)=0$.

If $v=1$, then $\partial D\cong P$ or $-P$.
The orientation on this boundary is the one induced by the
smooth spin filling $D$. Reversing orientation changes the
signature sign, so $\mu(-P)=-\mu(P)=1$ in $\Z/2$.
Theorem~\ref{input:rokhlin} would give
\[
 0=\sigma(D)\equiv8\mu(\pm P)=8\pmod{16},
\]
a contradiction.
Therefore $v=0$ and $\partial D\cong S^3$.
Consequently, $\pi_1(A)\cong\pi_1(\partial D)=1$, and
\[
 \pi_1(Y)\cong\pi_1(\partial D)*F_r\cong F_r.
\]
\end{proof}

\subsection{Proof of Lemma~\ref{lem:freeepi}}

\begin{proof}
Put $N=\ker a$.
Injectivity of $H_1(a;\Z)$ implies that every element of
$N$ has trivial image in the abelianization of $F_r$.
Thus
\[
 N\subset[F_r,F_r].
\]
Hopf's formula \cite[Section~0]{Brown} gives
\[
 H_2(G;\Z)
 \cong\frac{N\cap[F_r,F_r]}{[F_r,N]}
 =\frac{N}{[F_r,N]}.
\]
Since $H_2(G;\Z)=0$, it follows that
\[
 N=[F_r,N].
\]

Define the lower central series by
\[
 \gamma_1(F_r)=F_r,\qquad
 \gamma_{\ell+1}(F_r)=[F_r,\gamma_\ell(F_r)].
\]
Inductively,
\[
 N\subset\gamma_\ell(F_r)
 \quad\Longrightarrow\quad
 N=[F_r,N]\subset\gamma_{\ell+1}(F_r).
\]
\Needspace{7\baselineskip}
Hence
\[
 N\subset\bigcap_{\ell\ge1}\gamma_\ell(F_r)=1,
\]
where the last equality is residual nilpotence of free groups
\cite[Corollary~1.4]{MP}.
Thus $a$ is injective and therefore an isomorphism.
\end{proof}

\section{Proofs of Section~\ref{sec:assembly}}
\label{app:assembly}

\subsection{Proof of Proposition~\ref{prop:global-band}}

\begin{proof}
Apply Theorem~\ref{thm:cores} to $H=G$, and write the resulting
regular levels and cores as $t_i$ and $D_i$.
Their interiors exhaust $M$, so compactness of $\Gamma$
allows us to choose $j$ with
\[
 \Gamma\subset\intt D_j,\qquad t_j>0.
\]
Put $B=D_j$, $s=t_j$, and $a=G/s$.
Since
\[
 \Gamma=\{G\le0\},\qquad
 G|_{Y_\Gamma}=0,\qquad G|_{Y_B}=s,
\]
the function $a$ has the required boundary values.
Moreover, $0\le G\le s$ on $B\setminus\intt\Gamma$,
and regularity of the levels $0$ and $s$ gives
$da\ne0$ on $Y_\Gamma\cup Y_B$.

Positive rescaling preserves critical points, nondegeneracy,
and Morse indices. At a critical point $x$,
\[
 \Hess_q a=s^{-1}\Hess_qG,\qquad
 \ind_xa=\ind_xG\le2.
\]
At every regular point of the band,
\[
 \tr_{(\nabla_q a)^\perp}^q\Hess_q a
 =
 s^{-1}\tr_{(\nabla_q G)^\perp}^q\Hess_qG
 \ge s^{-1}\lambda_{(3)}^q(\Hess_qG)>0.
\]
Thus $a$ realizes $\mathsf{Band}_q(\Gamma,B)$.

Theorem~\ref{thm:cores} gives connected $Y_B$,
path-connected $O_B$, positive outward mean curvature,
and finite $\pi_1(O_B)$.
Theorem~\ref{thm:psc-defining} supplies the required
properties of $\Gamma$, including $H_q(Y_\Gamma)>0$.
Apply Proposition~\ref{prop:transfer} and
Corollary~\ref{cor:transferdata} with
\[
 H=G,\qquad K=\Gamma,\qquad D_j=B.
\]
They give a relative cell model for every prescribed finite
CW structure on $Y_\Gamma$ and show that $O_\Gamma$ is path connected,
with
\[
 |\pi_1(O_\Gamma)|<\infty.
\]
They also give
\[
 H_2(\Gamma)=H_3(\Gamma)=0,
\]
and show that $H_1(\Gamma)$ is torsion-free.
In particular, $\pi_1(O_\Gamma)$ contains no subgroup isomorphic
to $F_2$, and the relative model is available for the boundary
CW structure chosen in Proposition~\ref{input:normalform}.
The PSC boundary and inclusion-induced surjection
$\pi_1(Y_\Gamma)\twoheadrightarrow\pi_1(\Gamma)$
from Theorem~\ref{thm:psc-defining} supply the remaining
clauses of Condition~\ref{cond:topology}.
\end{proof}

\addtocontents{toc}{\protect\setcounter{tocdepth}{2}}

\end{document}